\documentclass[reqno,12pt]{preprint}

\usepackage[full]{textcomp}
\usepackage[osf]{newtxtext}
\usepackage{cabin}

\usepackage[varqu,varl]{inconsolata}
\usepackage[cal=boondoxo]{mathalfa}
\usepackage{microtype}

\usepackage{amsmath}
\usepackage{amssymb}
\usepackage{amsfonts}
\usepackage{hyperref}
\usepackage{breakurl}

\usepackage{mathrsfs}
\usepackage{enumitem}

\usepackage{tikz}
\usepackage{dsfont}
\usepackage{mhequ} 
\usepackage{mhsymb} 
\usepackage{mhenvs} 
\usepackage{array}

\usepackage{wasysym}
\usepackage{centernot}

\usetikzlibrary{snakes}
\usetikzlibrary{decorations}
\usetikzlibrary{positioning}
\usetikzlibrary{shapes}
\usetikzlibrary{external}

\colorlet{darkblue}{blue!90!black}

\DeclareFontFamily{U}{mathx}{\hyphenchar\font45}
\DeclareFontShape{U}{mathx}{m}{n}{
      <5> <6> <7> <8> <9> <10>
      <10.95> <12> <14.4> <17.28> <20.74> <24.88>
      mathx10
      }{}
\DeclareSymbolFont{mathx}{U}{mathx}{m}{n}
\DeclareMathSymbol{\bigtimes}{1}{mathx}{"91}

\def\emptyset{{\centernot\ocircle}}

\def\restr{\mathord{\upharpoonright}}
\def\id{\mathrm{id}}
\def\eqlaw{\stackrel{\text{\tiny law}}{=}}

\def\Sc{\mathcal{S}}

\colorlet{symbols}{blue!35!red}
\def\symbol#1{{\color{symbols}#1}}
\def\sX{\symbol{X}}
\def\sXi{\symbol{\Xi}}
\def\sone{\symbol{\mathbf{1}}}
\def\sPsi{\symbol{\Psi}}
\let\oldI\CI
\def\CI{\symbol{\oldI}}
\def\PPi{\boldsymbol{\Pi}}
\def\sscal#1{\langle\!\langle#1\rangle\!\rangle}
\def\ex{{\mathrm{ex}}}
\def\tree{\mathfrak{T}}
\def\dslash{/\kern-.37em/}

\newtheorem{example}[lemma]{Example}

\newtheorem{problem}{Open Problem}

\colorlet{testcolor}{green!60!black}

\tikzset{
	root/.style={circle,fill=testcolor,inner sep=0pt, minimum size=2mm},
	dot/.style={circle,fill=black,inner sep=0pt, minimum size=1mm},
	var/.style={circle,fill=black!10,draw=black,inner sep=0pt, minimum size=2mm},
	testfcn/.style={ultra thick,testcolor,shorten >=1pt,shorten <=1pt,<-},
	kernel/.style={semithick,shorten >=1pt,shorten <=1pt,->},
	kernel1/.style={semithick,shorten >=1pt,shorten <=1pt},
	kernel2/.style={->,semithick,shorten >=1pt,shorten <=1pt,postaction={decorate,decoration={markings,mark=at position 0.45 with {\draw[-] (0.05,-0.1) -- (0.05,0.1);\draw[-] (-0.05,-0.1) -- (-0.05,0.1);}}}},
	kernelBig/.style={semithick,shorten >=1pt,shorten <=1pt,decorate, decoration={zigzag,amplitude=1.5pt,segment length = 3pt,pre length=2pt,post length=2pt}},
	rho/.style={dotted,semithick,shorten >=1pt,shorten <=1pt},
	renorm/.style={shape=circle,fill=white,inner sep=1pt},
	labl/.style={shape=rectangle,fill=white,inner sep=1pt},
	xi/.style={circle,fill=symbols!10,draw=symbols,inner sep=0pt,minimum size=1.2mm},
	xix/.style={crosscircle,fill=symbols!10,draw=symbols,inner sep=0pt,minimum size=1.2mm},
	xib/.style={circle,fill=symbols!10,draw=symbols,inner sep=0pt,minimum size=1.6mm},
	xibx/.style={crosscircle,fill=symbols!10,draw=symbols,inner sep=0pt,minimum size=1.6mm},
	not/.style={circle,fill=symbols,draw=symbols,inner sep=0pt,minimum size=0.5mm},
	>=stealth,
	}
	
\tikzset{
	sdot/.style={circle,solid,draw=symbols,fill=symbols,inner sep=0pt,minimum size=0.5mm},
	yy/.style={circle,fill=gray!20,draw=black,inner sep=0pt,minimum size=0.8mm},
	}
	
\makeatletter
\def\DeclareSymbol#1#2#3{\expandafter\gdef\csname MH@symb@#1\endcsname{\tikz[baseline=#2,scale=0.15,draw=symbols]{#3}}%
\expandafter\gdef\csname MH@symb@#1s\endcsname{\scalebox{0.6}{\tikz[baseline=#2,scale=0.15,draw=symbols]{#3}}}}
\def\<#1>{\csname MH@symb@#1\endcsname}
\makeatother

\DeclareSymbol{X}{-2.4}{\node[sdot] {};}
\DeclareSymbol{1}{0}{\draw[white] (-.4,0) -- (.4,0); \draw (0,0)  -- (0,1.2) node[sdot] {};}
\DeclareSymbol{2}{0}{\draw (-0.5,1.2) node[sdot] {} -- (0,0) -- (0.5,1.2) node[sdot] {};}
\DeclareSymbol{3}{0}{\draw (0,0) -- (0,1.2) node[sdot] {}; \draw (-.7,1) node[sdot] {} -- (0,0) -- (.7,1) node[sdot] {};}
\DeclareSymbol{01}{-2}{\draw[white] (-.4,0) -- (.4,0); \draw (0,-1) -- (0.5,0)  -- (0,1) node[sdot] {};}
\DeclareSymbol{02}{-2}{\draw (-0.5,1) node[sdot] {} -- (0,0) -- (0.5,1) node[sdot] {}; \draw (0,0) -- (0,-1);}
\DeclareSymbol{03}{-2}{\draw (0,0) -- (0,1.2) node[sdot] {}; \draw (-.7,1) node[sdot] {} -- (0,0) -- (.7,1) node[sdot] {}; \draw (0,0) -- (0,-1);}
\DeclareSymbol{31}{-3}{\draw (0,0) -- (0,-1) -- (1,0) node[sdot] {}; \draw (0,0) -- (0,1.2) node[sdot] {}; \draw (-.7,1) node[sdot] {} -- (0,0) -- (.7,1) node[sdot] {};}
\DeclareSymbol{11}{-3}{\draw (.7,1) node[sdot] {} -- (0,0.3) -- (0,-1) -- (1,0) node[sdot] {};}
\DeclareSymbol{21}{-3}{\draw (0,0.3) -- (0,-1) -- (1,0) node[sdot] {}; \draw (-.7,1) node[sdot] {} -- (0,0.3) -- (.7,1) node[sdot] {};}
\DeclareSymbol{30}{-3}{\draw (0,0) -- (0,-1); \draw (0,0) -- (0,1.2) node[sdot] {}; \draw (-.7,1) node[sdot] {} -- (0,0) -- (.7,1) node[sdot] {};}
\DeclareSymbol{32}{-3}{\draw (0,0) -- (0,-1); \draw (1,0) node[sdot] {} -- (0,-1) -- (-1,0) node[sdot] {}; \draw (0,0) -- (0,1.2) node[sdot] {}; \draw (-.7,1) node[sdot] {} -- (0,0) -- (.7,1) node[sdot] {};}
\DeclareSymbol{22}{-3}{\draw (0,0.3) -- (0,-1); \draw (1,0) node[sdot] {} -- (0,-1) -- (-1,0) node[sdot] {};\draw (-.7,1) node[sdot] {} -- (0,0.3) -- (.7,1) node[sdot] {};}
\DeclareSymbol{20}{-3}{\draw (0,0) -- (0,-1);\draw (-.7,1) node[sdot] {} -- (0,0) -- (.7,1) node[sdot] {};}
\DeclareSymbol{32W}{-3}{\draw  [densely dotted, semithick] (-1.5,0)node[sdot] {}  -- (0,-1.5) node[yy]{} -- (1.5,0) node[sdot] {}; \draw (.8,-2.3)node{\tiny $y$} ;\draw[semithick] (0,0) -- (0,-1.5) node[yy] {}; \draw [densely dotted, semithick] (0,0) -- (0,1.8) node[sdot] {}; \draw[densely dotted, semithick] (-1,1.5) node[sdot] {} -- (0,0) -- (1,1.5) node[sdot] {};}
\DeclareSymbol{32WW}{-3}{\draw  [densely dotted, semithick] (-1.5,0)node[sdot] {}  -- (0,-1.5) node[yy]{} -- (1.5,0) node[sdot] {}; \draw (.8,-2.3)node{\tiny $\bar{y}$} ;\draw[semithick] (0,0) -- (0,-1.5) node[yy] {}; \draw [densely dotted, semithick] (0,0) -- (0,1.8) node[sdot] {}; \draw[densely dotted, semithick] (-1,1.5) node[sdot] {} -- (0,0) -- (1,1.5) node[sdot] {};}

\def\RR{\mathfrak{R}}
\def\TT{\mathscr{T}}

\def\simnot{\stackrel{\vbox to 0.15em{\hbox{\kern0.07em$^\circ$}}}{\sim}}

\begin{document}

\title{Renormalisation of parabolic stochastic PDEs \\[.4em] \small The 20th Takagi Lectures}
\author{M. Hairer}
\institute{Department of Mathematics, Imperial College London\\ \email{m.hairer@imperial.ac.uk}}

\maketitle

\begin{abstract}
We give a survey of recent result regarding scaling limits of 
systems from statistical mechanics, as well as the universality of the behaviour
of such systems in so-called cross-over regimes. It transpires that some of 
these universal objects are described by singular stochastic PDEs. We then give
a survey of the recently developed theory of regularity structures which allows to 
build these objects and to describe some of their properties. We place particular emphasis
on the renormalisation procedure required to give meaning to these equations.

These are expanded notes of the $20^{\mathrm{th}}$ Takagi lectures held at
Tokyo University on November 4, 2017.
\end{abstract}

\tableofcontents

\section{Scaling limits}

One major achievement of probability theory has been the construction of 
a number of universal objects that arise naturally as scaling limits of a
great number of natural probabilistic models. Let us start by going through
a number of examples in order to build some intuition on what scaling limits are, in
what context they arise, and what are their properties. 

\subsection{The Wiener process}

One of the simplest and most
ubiquitous such objects is the Wiener process, which is a random continuous
function $W \colon \R_+\to \R$ with the property that any two non-overlapping increments are
independent and that for any $t > s \ge 0$, $W(t) - W(s)$ is a centred Gaussian 
random variable with variance $t-s$. We further impose
that $W(0) = 0$.

The Wiener process arises naturally as a scaling limit in the following way.
Consider an i.i.d.\ sequence of random variables $\{X_n\}_{n \ge 0}$ that are
centred with finite variance $\sigma > 0$ and consider the partial sums
$S(n) = \sum_{k=0}^n X_k$ for $n \ge 0$. We extend $S$ 
to all arguments in $\R_+$ in a continuous way by linear interpolation: 
\begin{equ}[e:interpolate]
S(t) = \Bigl(\sum_{k=0}^{\lfloor t\rfloor} X_k\Bigr) + (t-\lfloor t\rfloor) X_{\lceil t\rceil}\;.
\end{equ}
We also introduce a family of rescaling operators
\begin{equ}
\bigl(\Sc_\lambda^{(\alpha)} f\bigr)(t) = \lambda^{-\alpha} f(\lambda t)\;.
\end{equ}
The following is then classical \cite{Donsker}:

\begin{theorem}\label{theo:CLT}
The family of random functions $S^\lambda = \Sc_\lambda^{(1/2)} S$
converges in law on $\CC(\R_+)$  to a Wiener process as $\lambda \to \infty$.
\end{theorem}

Observe that we start from a \textit{fixed} model (the random function $S$)
on which we perform a simple rescaling operation to obtain a family of models $S^\lambda$ depending
on a scale $\lambda$. This scale is then sent to infinity in order to obtain a scaling 
limit. It follows immediately that $W$ itself is invariant under the action of the
rescaling $\Sc_\lambda^{(1/2)}$ since
\begin{equ}
\Sc_\lambda^{(1/2)} W = \lim_{\mu \to \infty}\Sc_\lambda^{(1/2)}\Sc_\mu^{(1/2)} S
= \lim_{\mu \to \infty}\Sc_{\lambda+\mu}^{(1/2)} S
= W\;,
\end{equ}
where all limits and identities are in law.
The exponent $1/2$ appearing in this expression is called the
\textit{scaling exponent} of $W$. Another important feature of $W$ is that it exhibits a form
of \textit{stationarity}, or translation invariance. Indeed, the partial sums $S(n)$ 
are such that, for any fixed $k \ge 0$, the random sequence $\tilde S(n) = S(k+n) - S(k)$ is
equal in law to the random sequence $S$ itself. It immediately follows that $W$ satisfies the
analogous property, namely that $W(\cdot + s) - W(s)$ is equal in law to $W$ itself, for every $s \ge 0$.
A final and equally crucial property enjoyed by the Wiener process $W$ is the
\textit{Markov property}: for any time interval $I \subset \R_+$, the law of
$W \restr I$, conditional on $W \restr I^c$, only depends on the values of $W$ on
the boundary of $I$. This is again a simple consequence of the fact that $S$ satisfies a similar 
property, but with fattened boundaries. More precisely, for $[a,b] \subset \R_+$,
the law of $S \restr [a,b]$, conditional on $S \restr [a,b]^c$, only depends on
the values of $S$ on $[a-1,a] \cup [b,b+1]$.

\subsection{Critical percolation}

Scaling limits are of course not restricted to functions and, as a matter of fact,
several major recent results concern scaling limits of more sophisticated objects.
Let us briefly describe two of these results, namely critical percolation and the 
construction of the ``Brownian plane''. In the first example, our starting point is
the infinite regular triangular lattice, together with an assignment of a random
variable $\omega_x \in \{0,1\}$ for every vertex $x$ of the lattice. We assume that the $\omega_x$ are
i.i.d.\ Bernoulli random variables with parameter $p$. 

There are two different ways in which we can take a scaling limit in this case. The first ``na\"\i ve''
one, 
similarly to \eqref{e:interpolate}, is to extend $\omega$ to the whole plane by some
local interpolation procedure that respects the symmetries of the lattice. 
(Local in the sense that the value at a given point $x$ only depends on the values of 
$\omega_y$ for those vertices $y$ that intersect a neighbourhood of $x$ of some fixed radius.)
This yields a random function $\omega \colon \R^2 \to \R$ and we can ask whether this
random function admits a scaling limit. A variant of Theorem~\ref{theo:CLT} shows
that this is indeed the case, and we have the convergence in law
\begin{equ}[e:noise]
\lim_{\lambda \to \infty} \Sc_\lambda^{(-1)} (\omega - p)  = \xi\;,
\end{equ}
where $\xi$ is a ``white noise''. This time, the limit is not a random function,
but a random Schwartz distribution: for any collection of test functions
$\phi_k$, the random variables $\xi(\phi_k)$ are jointly centred Gaussian and
one has
\begin{equ}
\E \xi(\phi)\xi(\psi) = \scal{\phi,\psi}\;,
\end{equ}
where the scalar product is taken in $L^2$. In particular, if $\delta_\eps$ is an $\eps$-approximation
of a Dirac mass, one has $\E \xi(\delta_\eps)^2 \sim \eps^{-2}$, which shows that there is no hope
to be able to find a continuous version of $\xi$.
The fact that white noise can only be realised as a random distribution is directly linked to the
fact that it arises as a scaling limit \eqref{e:noise} with a \textit{negative} exponent.

However, while this white noise scaling limit does describe the large-scale behaviour of
local averages of the random variables $\omega_x$, it does
not tell us anything about the large-scale behaviour of many interesting observables that are
naturally built from $\omega$. In particular, one would like to be able to give a description
of the large-scale geometry of the subset $\{x\,:\, w_x = 1\}$ of the triangular lattice, together 
with the connectivity structure induced by that of the underlying lattice. 
One possible way of encoding this connectivity structure is by viewing it as a collection of
``quads'' (i.e.\ homeomorphisms $Q\colon [0,1]^2 \to \R^2$), where we keep precisely those quads 
$Q$ such that there exists a percolation cluster intersecting the two opposite edges
$Q(\{0\} \times [0,1])$ and $Q(\{1\} \times [0,1])$. (We can view the percolation clusters as closed
subsets of the plane by for example linking neighbouring open sites with a line segment. The precise
way in which this is done is irrelevant for what follows.)
Writing $\CQ$ for the space of all quads, which can be endowed with a natural distance function,
a percolation configuration is therefore encoded by a random subset $S_\omega \subset \CQ$
satisfying furthermore a number of natural consistency and monotonicity properties. 
The space $\CX$ of all such subsets can itself be endowed with a natural topology.
The is furthermore a natural scaling operation $\Sc_\lambda \colon \CQ \to \CQ$
on the space of quads by setting
\begin{equ}
\bigl(\Sc_\lambda Q\bigr)(x) = \lambda\, Q(x)\;. 
\end{equ}
This in yields a map $\Sc_\lambda \colon \CX \to \CX$ on subsets of $\CQ$ which preserves
the consistence and monotonicity properties alluded to above.
It was then shown by Schramm and Smirnov \cite{SSPerco} that if one sets $p = p_c = {1\over 2}$ (the
critical value for percolation on the triangular lattice \cite{Kesten}), then 
the laws of the sequence of random sets $\Sc_\lambda S_\omega$ is tight  in $\CX$
as $\lambda \to \infty$ and its accumulation points are non-degenerate. 
The limiting $\CX$-valued random variable  $\tilde S$ (which is conjectured to be unique) 
has again natural scale invariance,
 translation invariance, and Markov properties. In this particular case, it is furthermore
rotation invariant and actually conformal invariant. 

It is interesting to note that the information encoded in the scaling limit $\xi$ and that encoded
in the scaling limit $\tilde S$ are completely different, although both $\Sc_\lambda^{(-1)} w$
and $\Sc_\lambda S_\omega$ encode the exact same information
for any fixed $\lambda > 0$. What ``completely different'' means in this context is
that the noise (in the sense of Tsirelson \cite{Tsirelson}) generated by $\xi$ is ``white'' while
that generated by $\tilde S$ is ``black'' \cite{SSPerco}.

\subsection{The Ising model}

Another example of scaling limit was recently obtained by 
\cite{Ising1,Ising2} for the 2D Ising model at criticality. Recall that the Ising model in a
domain $\Lambda \subset \Z^2$ is given by the measure on $\{\pm 1\}^\Lambda$
assigning to each configuration $\sigma$ a probability proportional to 
$\exp(-\beta H(\sigma))$, where the Hamiltonian is given by 
\begin{equ}
H(\sigma) = - \sum_{x \sim y} \sigma_x \sigma_y\;.
\end{equ}
Here, $x\sim y$ if and only if $x$ and $y$ are nearest neighbours in $\Z^2$ and
the values of $\sigma$ outside of $\Lambda$ are considered to be fixed.
It is well-known that the Ising model exhibits a phase transition: there exists a
value $\beta_c > 0$ such that, for $\beta \le \beta_c$, there exists a \textit{unique}
probability measure on $\{\pm 1\}^{\Z^2}$ with the property that, for every
$\Lambda \subset \Z^2$ finite, the conditional measure on $\{\pm 1\}^\Lambda$
is as above. For $\beta > \beta_c$ however, while such measure do still exist, they
are not unique anymore. In fact, there are \textit{two} such measures, one
in which a majority of spins take the value $+1$, and one in which a majority of spins take the
value $-1$. At criticality (i.e.\ for $\beta = \beta_c$), local averages of spins exhibit 
a non-trivial scaling limit in the sense that if we extend $\sigma$ to all of $\R^2$ in a
way similar to above, then $\Sc_\lambda^{(-1/8)} \sigma$ converges in law to a non-trivial
measure $\P_c$ on the space of distributions on $\R^2$. Unlike in the previous
examples of this type, the measure $\P_c$ is not Gaussian, in fact its tails are lighter
than Gaussian.
Again, the measure $\P_c$ is scaling invariant (with exponent $-1/8$), translation invariant,
and satisfies the spatial Markov property.

\subsection{Interface fluctuations}
\label{sec:interface}

Our final example of scaling limit is the closest one to the type of 
problems considered in these notes, and this is also the only example of scaling limit
in which one considers models that depend on both space and time.
Consider a model of one-dimensional interface growth where the interface is modelled 
as the graph of a function $h\colon \Z \to \Z$. We furthermore restrict the state space to those
functions $h$ such that $h(x+1) - h(x) \in \{\pm 1\}$ for every $x \in \Z$.
A simple dynamic on such functions is given by the following. To each site $x \in \Z$,
we then assign a Poisson process $t \mapsto N_x(t)$, with all the $N_x$'s being i.i.d.\ and having the same rate
(say $1$). This allows us to build a dynamic on the space of height functions $h$ in the
following way. Whenever one of the processes $N_x$ jumps, we update $h$ in the following way.
If $h_t$ has a local minimum at $x$, then we set $h_{t^+}(x) = h_t(x) + 2$.
If $h_t$ has a local maximum at $x$, then we set $h_{t^+}(x) = h_t(x) - 2$, but only with
some fixed probability $q \in [0,1]$, independently of the processes $N_x$.
If $h_t$ has neither a local minimum, nor a local maximum at $x$, then we leave it unchanged.
Note that the case $q=1$ is special since in this case $h$ and $-h$ have the same distribution. 

It is possible to show \cite{Liggett} that this dynamic is indeed well-defined
(this is not completely obvious since infinitely many of the processes $N_x$ perform
a jump in any given time interval). Depending on the literature, this process is referred to as
the height function for the asymmetric simple exclusion process (symmetric in the case $q=1$)
or as the SOS (Solid On Solid) model.
This process has the property that, 
for any given value of $q$, the simple random walk is an invariant measure for it.
This suggests that, in order to obtain a scaling limit, the initial condition $h_0$ should
be chosen such that $\Sc_\lambda^{(1/2)} h_0$ converges in law to some limit $\tilde h_0$.
In this situation, since time and space play different roles, there is no reason a priori 
to consider time scales of the same order as the spatial scale, so we consider
scaling operators with two scaling exponents like so:
\begin{equ}
\bigl(\Sc_\lambda^{(\alpha,\beta)}h\bigr)(x,t) = \lambda^\alpha h(\lambda x, \lambda^\beta t)\;.
\end{equ}
In the symmetric case $q = 1$, one then has the following result \cite{SSEP1,SSEP2}:

\begin{theorem}
For $q=1$, let $h_0^\lambda$ be a sequence of initial conditions for the symmetric SOS model such that 
$\Sc_\lambda^{(1/2)} h_0^\lambda$ converges to a limit $h_0$ in $\CC_b(\R)$. Then, the process $\Sc_\lambda^{(1/2,2)}h$
converges in law as $\lambda \to \infty$ to the solution to the stochastic heat equation with initial condition $h_0$:
\begin{equ}[e:SHE]
\d_t h = \d_x^2 h + \xi\;,
\end{equ}
where $\xi$ denotes space-time white noise.
\end{theorem}

Note again that the process $h$ satisfies the (space-time) Markov property as a consequence of the fact
that it has a purely local specification. In the case of interface fluctuation models, it turns out that 
symmetric models behave very differently at large scales from asymmetric ones.
Indeed, the following has very recently been shown in \cite{KPZ}.

\begin{theorem}
For $q=0$, let $h_0^\lambda$ be a sequence of initial conditions for the totally asymmetric SOS model such that 
$\Sc_\lambda^{(1/2)} h_0^\lambda$ converges to a limit $h_0$ in $\CC_b(\R)$. Then, the process $\Sc_\lambda^{(1/2,3/2)}h$
converges in law as $\lambda \to \infty$ to a limiting non-trivial Markov process $h$.
\end{theorem}

Note that the time-scaling exponent appearing here is different from the one appearing 
in the symmetric case. Furthermore, \cite{KPZ} give a full characterisation of the transition
probabilities of the limiting process $h$, which shows in particular that it is not Gaussian.
Unlike in the symmetric case, the scaling limit in the totally asymmetric case (also called
the ``KPZ fixed point'') is not described as the solution to a stochastic PDE. It does however
have a characterisation via a type of noisy Lax-Oleinik formula \cite{AirySheet}, although this characterisation
is not complete. (There are only tightness results and partial characterisations of the ``Airy sheet''
which acts as a driving noise.)

\subsection{General features and open problem}

Despite the impression the reader may obtain from these various 
examples, these scaling limits are rigorously much less understood than one would think.
For example, the Ising model admits two very natural dynamics preserving its Gibbs measure: 
the Glauber dynamic in which one randomly flips spins and then rejects some moves in order
to enforce the correct invariant measure and the Kawasaki dynamic, where the elementary moves are given
by exchanging neighbouring spins. It is then natural to ask whether there is a scaling limit
for the process $(t,x) \mapsto \sigma_x(t)$ induced by these dynamics at the critical 
parameter $\beta_c$. It is conjectured that this is indeed the case, but one does not even
have a conjecture for the correct value of the corresponding dynamical scaling exponents! 
Similarly, one does not know what the correct scaling $\alpha$ is which guarantees that
$\Sc_\lambda^{(-\alpha)} \sigma$ converges to a non-trivial limit for $\sigma$ given
by the 3D Ising model at criticality, although recent breakthrough results \cite{Bootstrap1,Bootstrap2} 
yield a constructive, albeit non-rigorous,
way of approximating this exponent (the best known approximation to date is 
about $0.51815$ which does not appear to be a rational number with small denominator).

Another feature of these scaling limits that still lacks rigorous mathematical understanding is
their universality. In all of the examples given above, we have mathematical theorems exhibiting
scaling limits for one specific model, typically one that is exactly solvable. 
It is (loosely) conjectured that in all cases, the \textit{same} scaling limit arises when 
starting from any model that exhibits the same broad features as the specific model 
we consider. This is only well understood in the first example of the Wiener process, where the functional central
limit theorem was shown to hold for a very large class of models.

In the case of critical percolation for example, one expects to obtain the same
scaling limit for its connectivity structure, independently of the precise structure of the
underlying graph, provided that it is planar and sufficiently regular. In the case of the
Ising model, one also expects the same results for a large class of underlying grids,
similarly to the case of percolation.
(For a different class of observables, this was indeed shown in \cite{UniversalIsing}, although the 
class of grids considered there is still rather rigid and interactions are still restricted to nearest-neighbours
in order to enforce a suitable type of integrability.) 
However, when considering the scaling limit of the magnetisation
as described above, one also expects many continuous models to fall into the same
universality class. For example, the $\Phi^4_2$ measure constructed in \cite{Nelson,GlimmJaffe,BrydgesAl}
is a measure on distributions on $\R^2$ which can formally be thought of as having 
density proportional to 
\begin{equ}
\exp \Bigl(- \int_{\R^2} \bigl(|\nabla \Phi(x)|^2 + \Phi^4(x) - C \Phi^2(x)\bigr)\,dx \Bigr)\;,
\end{equ}
with respect to ``Lebesgue measure'' (which of course does not exist).
This models the same situation as that of the Ising model and depends on a parameter 
$C$ which plays a role similar to that of the inverse temperature $\beta$. It is natural to conjecture that,
when this parameter is set equal to its critical value, the large-scale behaviour of the $\Phi^4_2$ measure
is precisely the same as that of the Ising model.

Finally, in the case of interface fluctuations, one expects the KPZ fixed point to 
be the scaling limit of a very large class of asymmetric models. In particular, 
one expects it to arise as the scaling limit of the SOS model for any $q \in [0,1)$,
but also as the scaling limit of various other models like the Eden model, the 
ballistic deposition model, etc. The only rigorous results in this direction were
obtained for models exhibiting a type of ``integrable structure'', for which it is possible
to show that certain observables converge to the corresponding observables of the KPZ
fixed point. See for example \cite{Kurt,TW,Macdonald} for some results in this direction.

All this suggests that answering the following open problem would allow to make substantial progress.

\begin{problem}
For any $d \ge 2$, 
characterise all pairs $(\alpha, \eta)$ such that $\alpha \in \R$ and $\eta$ is a
random distribution on $\R^d$ which is stationary, satisfies $\Sc_\lambda^{(-\alpha)} \eta = \eta$ in
law, and satisfies the spatial Markov property.
\end{problem}

\begin{remark}
If $\alpha \ge 0$, one cannot in general expect $\eta$ to be both self-similar with exponent $\alpha$ 
and stationary. However, at least for non-integer values, one can still expect stationarity in a generalised 
sense where one recenters $\eta$ in such a way that its value and the relevant number of derivatives vanish
at the origin. This is precisely what happens in the case of the two-sided Wiener process for example 
which is invariant under the action of the group of translations given by
\begin{equ}
(\tau_s W)(t) = W(t+s) - W(s)\;.
\end{equ}
\end{remark}

\begin{remark}
There are many natural variants of this question. For example, it would be natural
to also impose rotation invariance. It would also be natural to add time as an additional
dimension and to allow the scaling exponent to be different for time
scalings. Finally, it would be natural to impose restrictions on the tail behaviour of $\eta$,
for example that $\eta(\phi)$ has moments of all orders for every
$\phi \in \CC_0^\infty$. (In the case $d=1$, this property is known to significantly reduce 
the possible candidates.)
\end{remark}

\section{Crossover regimes}

We have seen in the previous section that there is a great interest in finding and characterising 
scaling limits, as well as the corresponding universality classes. (The universality class for a given
scaling limit $\eta$ consists of the collection of all those models that converge to $\eta$ under
a suitable rescaling operation which leaves $\eta$ invariant.)
Unfortunately, when considering scaling limits of ``random field'' type, the only 
situation that is mathematically well-understood in the sense that we have a mathematically rigorous 
understanding of
both the scaling limit itself and of its universality class is the Gaussian case. 
There are some non-Gaussian cases that can be characterised by explicit formulae, which 
include the case of two-dimensional models 
(without time-dependence), as well as the KPZ fixed point. In these cases however, universality
statements are restricted to microscopic models that share some of the integrability properties of 
their scaling limits.

This suggests that one possible way of improving our understanding of the large-scale behaviour
of random systems is to consider what happens ``in the vicinity'' of systems that rescale to
a Gaussian scaling limit. More precisely, consider a family $\{S^{(\eps)}\}_{\eps \in D}$
of random functions on $\R^d$ depending continuously on $\eps$, where $D \subset \R^m$ is some 
closed set containing the origin (think of $D = [0,1]$). Assume that this family is 
such that, for some scaling exponent $\alpha$, 
$\CS_\lambda^{(\alpha)} S^{(0)}$ converges, as $\lambda \to \infty$, to a Gaussian scaling limit $\eta$.
One can then find situations in which this is no longer the case for $\eps \neq 0$, so that 
one \textit{expects} to have some different scaling exponent $\beta \neq \alpha$ (or pair of 
space-time scaling exponents) such that, for fixed
$\eps \neq 0$, $\CS_\lambda^{(\beta)} S^{(\eps)}$ converges, as $\lambda \to \infty$, to some
non-Gaussian scaling limit. (It is usually completely out of reach of current mathematical technology
to prove such a thing. In general, it may also happen that $D$ is further broken into subregions on which 
$\beta$ takes different values.) It is then natural to ask whether it is possible to prove 
a non-trivial limiting result for limits of the form $\CS_{\lambda_\eps}^{(\alpha)} S^{(\eps)}$
as $\eps \to 0$ with $\lim_{\eps \to 0} \lambda_\eps = \infty$.

By a simple diagonal argument, using the continuity of $\eps \mapsto S^{(\eps)}$,
if $\lambda_\eps$ diverges sufficiently slowly, one still has $\lim_{\eps \to 0}
\CS_{\lambda_\eps}^{(\alpha)} S^{(\eps)} = \eta$. On the other hand, since $\beta \neq \alpha$,
one expects $\CS_{\lambda_\eps}^{(\alpha)} S^{(\eps)}$ to either converge to $0$ or blow up 
as $\eps \to 0$ for $\lambda_\eps$ diverging sufficiently fast. This suggests 
the existence of a non-trivial intermediate regime, called the \textit{crossover regime} 
for which some non-trivial limits different from $\eta$ are obtained. 

In these notes, we'll argue that in many cases of interest, the behaviour one observes is
the following. There exists a family $\eta^v$ of possible limits, parametrised by $v$
belonging to some closed subset $K \subset \R^m$ containing the origin for some $m \ge 0$, 
with $\eta^0 = \eta$, the original
Gaussian scaling limit. Each element of the family is typically stationary (possibly in a 
generalised sense as already described earlier) and  
Markovian, but it is this time only the family as a whole which is scaling invariant
and not individual members of the family.
More precisely, there exists an action $T$ of $(\R_+,\cdot)$ onto $K$
such that
\begin{equ}[e:action]
\CS_{\lambda}^{(\alpha)} \eta^v = \eta^{T^\lambda v}\;,
\end{equ}
and such that $\lim_{\lambda \to 0} T^\lambda v = 0$ for every $v \in K$.
In other words, the family is invariant under the scaling operations with exponent $\alpha$
and, at \textit{small} scales, every member of the family looks like $\eta$.
Furthermore, the family $\{\eta^v\}_{v \in K}$ is \textit{universal} in the sense that, for a
very large class of models $\{S^{(\eps)}\}_{\eps \in D}$ as above, one can find a map
$\eps \mapsto (\lambda_\eps, v_\eps)$ such that 
\begin{enumerate}[itemsep=0pt,parsep=0pt,topsep=0.2em]
\item $v_\eps$ is uniformly bounded
on $\eps \in D \setminus \{0\}$,
\item there are sequences $\eps_n \to 0$ such that $v_{\eps_n}$ remains uniformly bounded away from $0$, 
\item $\lambda_\eps \to \infty$ as $\eps \to 0$,
\item one has $\CS_{\lambda_\eps}^{(\alpha)} S^{(\eps)} \to \eta^{v_\eps}$ as $\eps \to 0$.
\end{enumerate}
Note that there is necessarily some indeterminacy here since, if $\eps \mapsto (\lambda_\eps, v_\eps)$
is one map satisfying the above properties then it follows from \eqref{e:action} that, 
for any fixed $\mu > 0$, 
the map $\eps \mapsto (\mu \lambda_\eps, T^\mu v_\eps)$ also satisfies the same properties.

\begin{remark}
In all the examples below, the underlying space $\R^d$ on which our random distributions are
defined naturally breaks into a time component and $d-1$ spatial components. In this case, scalings
are typically different for these two types of components, but the discussion above remains
unchanged otherwise.
\end{remark}

\begin{remark}
There are some situations where major technical difficulties arise due to the fact that 
space is unbounded. Although this is very much case-dependent, there are often natural 
ways of defining the models $S^{(\eps)}$ and $\eta^v$ on a finite domain rather than the 
whole space. In this case, the same discussion usually applies but with two caveats.
First, there may be additional parameters governing boundary conditions and matching these between
$S^{(\eps)}$ and $\eta^v$ may be non-trivial, see \cite{CS,Mate}. Second, 
if we normalise the $\eta^v$ to be defined on a domain of size $1$, then the models
$S^{(\eps)}$ have to be defined on a domain of size $\lambda_\eps$ for the third property above to
hold.
\end{remark}

We now discuss a couple of concrete examples to illustrate the situation.

\subsection{The KPZ equation}

Let us now show how the KPZ equation \cite{KPZOrig} fits into the
framework described above. Recall that the KPZ equation is the stochastic partial differential equation
\textit{formally} given by
\begin{equ}
\d_t h = \d_x^2 h + (\d_x h)^2 + \xi\;,
\end{equ}
where $\xi$ denotes space-time white noise. The reason why this expression is only formal is that 
solutions are not differentiable, so the expression $(\d_x h)^2$ has no classical meaning.
It is however possible to give an unambiguous notion of ``solution'' to this equation,
the so-called Hopf-Cole solution, see \cite{BG,JeremyHerbert}. As a matter of fact, this gives a meaning
the family of equations formally given by
\begin{equ}[e:KPZ]
\d_t h = \d_x^2 h + c_1(\d_x h)^2 + c_2 + \xi\;,
\end{equ}
and we write $\eta^c$ with $c = (c_1,c_2) \in \R^2$ for its solutions. 

\begin{remark}
There is a rather non-trivial question arising here which is that of the choice of
initial condition. It so happens that it is possible to construct solutions
that are stationary in space-time, modulo recentering, and these are the solutions that we consider.
There is a subtlety here in that these solutions are themselves not unique since
Brownian motion with drift is invariant for the corresponding Markov process for any value of the drift,
see \cite{FQ}. We lift this ambiguity by choosing the (unique) solution with zero drift.
\end{remark}

Using the scale invariance properties of white noise,
it can be verified that, at least at a formal level, the family of 
solutions to \eqref{e:KPZ} does satisfy \eqref{e:action} with 
\begin{equ}[e:actionKPZ]
T^\lambda(c_1,c_2) = (\lambda^{1/2} c_1, \lambda^{3/2} c_2)\;.
\end{equ}
It happens that the Hopf-Cole solutions are such that this identity does indeed hold
for them as well.

\begin{remark}
We could of course also add constants in front of the terms $\d_x^2 h$ and $\xi$, but it
can be shown that it is always possible to rescale solutions in such a way that these constants
are equal to $1$. As a matter of fact, even the constants $c_1$ and $c_2$ can be eliminated by
performing first a rescaling to set $c_1$ to $1$ by \eqref{e:actionKPZ} and then perform
a height shift, i.e.\ considering transformations of the type
$h(t,x) \mapsto h(t,x) - c t$, which allow to set $c_2$ to $0$ without loss of generality.
This is why we usually talk about ``the'' KPZ equation, rather than the family of KPZ equations.
\end{remark}

\subsubsection{Convergence of the SOS model}

Let now $S^{(\eps)}$ denote the stationary SOS model with $q = 1-\sqrt \eps$ and with
fixed time distribution given by that of the simple random walk. We furthermore perform
a height shift in order to set the average speed of the interface to zero, so that if
$h^{(\eps)}$ denotes the stationary SOS model with $q = 1-\sqrt \eps$ as defined in
Section~\ref{sec:interface}, we set
\begin{equ}
S^{(\eps)}(t,x) = h^{(\eps)}_t(x) - \sqrt \eps t\;.
\end{equ}
It was then shown in \cite{BG} that if we choose $\lambda_\eps = 1/\eps$, 
and consider the processes $h^{(\eps)} = \CS_{\lambda_\eps}^{(1/2,2)}S^{(\eps)}$,  rescaled as in 
the \textit{symmetric} case mentioned in Section~\ref{sec:interface}, then 
$h^{(\eps)}$ does converge to the Hopf-Cole solution to the KPZ equation. 

This can indeed be cast into the framework exposed at the beginning of this section: 
We have $D = [0,1]$ and $K = \R^2$, and the mapping $\eps \mapsto (\lambda_\eps, v_\eps)$ is 
simply given by $\eps \mapsto (\eps^{-1}, c)$, where $c$ denotes the specific value of the 
constants in the KPZ equation appearing in the result of \cite{BG}.
In this case, this seems very much overkill. The interesting fact however is that 
the limiting object obtained in the crossover regime, namely the solution to the KPZ
equation, is itself \textit{universal} in the sense that it arise in many more situations. 

\subsubsection{A universality result}

Consider for example a Gaussian random field $\zeta\colon \R^2 \to \R$ with a correlation
function $\rho$ that is compactly supported and such that $\int \rho = 1$. Then, a simple
calculation shows that $\CS^{(-3/2,2)}_\lambda \zeta$ converges to space-time white noise
as $\lambda \to \infty$. In particular, if we denote by $h^{(0)}$ the solution to
\begin{equ}[e:basic]
\d_t h^{(0)} = \d_x^2 h^{(0)} + \zeta\;,
\end{equ}
then $\CS^{(1/2,2)}_\lambda h^{(0)}$ converges in law as $\lambda \to \infty$ to the solution
to the stochastic heat equation \eqref{e:SHE}. Consider now perturbations of \eqref{e:basic}
of the type
\begin{equ}[e:interface]
\d_t h = \d_x^2 h + P(\d_x h) + \zeta\;,
\end{equ}
where $P$ is an even polynomial of some fixed degree $2m$. 
The rationale for restricting ourselves to even polynomials is that it is natural
for interface growth models to be symmetric under the exchange $x \leftrightarrow -x$.
Furthermore, the rationale for considering perturbations depending only on $\d_x h$ 
and not on $h$ itself is that the environment in which our interface moves is 
homogeneous.

This again fits the framework described above, but this time it is the polynomial
$P$ itself which plays the role of the parameter $\eps$, so that one has
$D \approx \R^{m+1}$, with the identification between polynomials and $\R^{m+1}$ given by
\begin{equ}
P = P_\eps\;,\qquad P_\eps(u) = \sum_{p=0}^m \eps_p u^{2p}\;.
\end{equ}
What behaviour does one expect for the large-scale behaviour of \eqref{e:interface}
for small $P$? Note that while \eqref{e:basic} only admits a stationary solution in a
generalised sense, its spatial derivative $\d_x h^{(0)}$ admits a genuine stationary
solution. Since it is Gaussian, this means that there exists a centred Gaussian measure
$\nu_0$ such that $\d_x h^{(0)}(t,x)$ has law $\nu$ for any fixed $(t,x)$. 

It is then not unreasonable to expect that the effective behaviour of $P(\d_x h)$ in 
\eqref{e:interface} is close to that of the constant $\int P(u)\,\nu(du)$.
Indeed, it is not difficult to show that this is the case in the sense that if
we choose for example
\begin{equ}[e:lambdaeps]
\lambda_\eps = |P_\eps|^{-2/3}\;,
\end{equ}
where $|P_\eps| = \sum_p |\eps_p|$, then $\CS^{(1/2,2)}_{\lambda_\eps} h^{(P_\eps)}$
is uniformly close for $|P_\eps|$ small to the solution to
\begin{equ}[e:SHEshift]
\d_t h = \d_x^2 h + c_\eps + \xi\;, \qquad c_\eps \eqdef \lambda_\eps^{3/2} \int P(u)\,\nu(du)\;.
\end{equ}
Here, we wrote $h^{(P)}$ for the solution to \eqref{e:interface}.
Since \eqref{e:lambdaeps} guarantees that $c_\eps$ is uniformly bounded, this is indeed
a result of the type described above, but it is somewhat disappointing. Indeed, as already mentioned,
the solution to \eqref{e:SHEshift} differs from that of \eqref{e:SHE} by a simple shift,
so that we really haven't obtained a different object at all. 

Instead, we should focus on what happens when $\int P\,d\nu$ vanishes, in which
case $c_\eps$ vanishes, so that we have a chance of being able to look at larger scales.
It turns out that a more interesting choice of scale $\lambda_\eps$ which leads to non-trivial
behaviour is given by 
\begin{equ}[e:betterscale]
\lambda_\eps^{-1} = |P_\eps|^{2} + \Bigl|\int P(u)\,\nu(du)\Bigr|^{2/3}\;.
\end{equ}
We see that if $\int P\,d\nu$ is of the same order as $|P_\eps|$, then this scale is 
of the same order as that in \eqref{e:lambdaeps}. If however $\int P\,d\nu$ is much smaller,
then this scale is much larger than the previous one.

\begin{remark}
The precise form of $\lambda_\eps$ does not matter as long as it satisfies the bounds
$\lambda_\eps \lesssim |P_\eps|^{-2}$ and $\lambda_\eps \lesssim |\int P\,d\nu|^{-2/3}$.
\end{remark}

The main result of \cite{Jeremy} can then be loosely formulated as follows.

\begin{theorem}\label{theo:univKPZ}
There exists a trilinear form $F\colon P \mapsto F(P,P,P)$ such that,
with the choice \eqref{e:betterscale}, $\CS^{(1/2,2)}_{\lambda_\eps} h^{(P_\eps)}$
is uniformly close for $|P_\eps|$ small to the Hopf-Cole solution to \eqref{e:KPZ}
with
\begin{equ}
c_1 = {\lambda_\eps^{1/2} \over 2} \int P''(u)\,\nu(du)\;,
\qquad
c_2 = \lambda_\eps^{3/2} \Bigl(\int P(u)\,\nu(du) +  F(P,P,P)\Bigr)\;.
\end{equ}
\end{theorem}

Similar results have recently been obtained in \cite{MaxNikolas}, where the authors also allow
for non-polynomial $P$, but with the caveat that they are restricted to a slightly different class
of models for which the invariant measure is explicitly known. The work \cite{Weijun} 
also extends the above results to general $P$, this time for the same model as \eqref{e:interface}.

\begin{remark}
The formulation of Theorem~\ref{theo:univKPZ} is not complete as it leaves out the question of the initial condition.
Furthermore, the results just mentioned only hold for \eqref{e:interface} defined on a torus of size
$\lambda_\eps$ instead of the whole real line. The reason is that, on the whole real line, it is not even clear
for which initial conditions \eqref{e:interface} admits local solutions!
\end{remark}

These results suggest the following conjecture which would be a very strong universality 
statement for the KPZ equation.

\begin{problem}
Show that the solution to the KPZ equation is the {\em unique} continuous
stationary (in the generalised sense) space-time Markov process $h$ invariant under $x \leftrightarrow -x$ 
which has the property that 
$\CS^{(1/2,2)}_\lambda h$ converges to the solution to the stochastic heat equation
as $\lambda \to 0$ and $\CS^{(1/2,3/2)}_\lambda h$ converges to the KPZ fixed point
as $\lambda \to \infty$.
\end{problem}

\subsection{The dynamical \texorpdfstring{$\Phi^4_3$}{Phi43} model}

This is a family of models which should formally be thought of as the solutions
to the stochastic PDEs
\begin{equ}[e:Phi4]
\d_t \Phi = \Delta \Phi + c_1 \Phi - c_2 \Phi^3 + \xi\;,
\end{equ}
where $\xi$ again denotes space-time white noise, but the spatial variable now takes 
values in $\R^3$. This time, the model with $c = (c_1,c_2) =0$ is invariant under the action of
the scaling operators $\CS^{(-1/2,2)}_{\lambda_\eps}$.
A formal calculation suggests that this family of models satisfies \eqref{e:action}
with the action $T^\lambda$ given by 
$T^\lambda c = (\lambda^{3/2} c_1, \lambda^{1/2} c_2)$. This is \textit{incorrect} however,
instead it is the case that it does satisfy \eqref{e:action}, but with the modified 
action
\begin{equ}
T^\lambda c = \big(\lambda^{3/2} (c_1 +  a c_2^3 \log \lambda), \lambda^{1/2} c_2\big)\;,
\end{equ}
where $a$ is some fixed explicitly computable constant. 
(One readily verifies that this action satisfies the semigroup
property $T^\lambda T^\mu = T^{\lambda \mu}$.)

This time, while the constant $c_2$ can always be set to $1$ by a simple rescaling operation,
there is no way in general to adjust the value of $c_1$ by a simple local transformation,
so we do genuinely have a one-parameter family of distinct models indexed by $c_2$.
This parameter plays a role similar to that of the temperature in the Ising model, and the
two models are expected to fall into the same universality class in the sense that one expects
there to exist a critical value for $c_2$ at which \eqref{e:Phi4} admits a scaling limit which coincides
with the scaling limit for the Glauber dynamic of the critical 3D Ising model.
Let us reiterate again that although this is an accepted fact in the physics
literature, from a mathematical perspective it is highly speculative. In particular, neither of 
these scaling limits is known to exist, and there does not even exist a precise conjecture 
for the corresponding scaling exponents.

It was shown in \cite{WeijunPhi4} that, if one considers models of the type
\begin{equ}[e:potential]
\d_t \Phi = \Delta \Phi - V'(\Phi) + \zeta\;,
\end{equ}
for even polynomials $V\colon \R \to \R$ of fixed degree and $\zeta$ a Gaussian random field
with compactly supported correlations as above, then a result analogous to \eqref{theo:univKPZ}
holds (but with slightly different expressions for $c$ and the scale $\lambda$).
It was also shown in \cite{HendrikJC} that the two-dimensional analogue of this model
arises in the crossover regime for the Glauber dynamic of Ising-Kac models, namely
spin models with Hamiltonian
\begin{equ}
H_\eps(\sigma) = - \eps^2 \sum_{x,y} \sigma_x \sigma_y K(\eps(x-y))\;.
\end{equ}
A similar result is expected to hold in the three-dimensional case as well.

\section{Singular stochastic PDEs}

In the previous section, we have seen singular stochastic PDEs appearing as universal objects
describing the crossover regimes for various models from statistical mechanics. 
These equations were singular in the sense that they involved nonlinearities that appear to have no
canonical meaning since they involve products of distributions of negative order.
For this reason, it is not clear how processes like the solutions to (\ref{e:KPZ},\ref{e:Phi4}) 
should be defined and, a fortiori, how universality results like the ones presented
above should be obtained.

\subsection{Analysis of the \texorpdfstring{$\Phi^4$}{Phi4} model}
\label{sec:Phi4}

Let us focus now on \eqref{e:Phi4} with $c_2 = 1$ since the KPZ equation has non-generic
additional structure which makes it amenable to alternative techniques \cite{KPZEnergy,KPZEnergyUnique}.
We are therefore interested in the study of the family of solutions to
\begin{equ}[e:Phi4bis]
\d_t \Phi = \Delta \Phi + c \Phi - \Phi^3 + \xi\;,
\end{equ}
with $c \in \R$. The problem is that solutions to the stochastic heat equation
\begin{equ}[e:SHEadd]
\d_t X = \Delta X + \xi\;,
\end{equ}
in three dimensions (and therefore presumably also solutions to \eqref{e:Phi4bis}, whatever meaning
we choose to give to the equation) are random distributions of order $-1/2$.
In other words, the operation $\Phi \mapsto \Phi^3$ is not closable on any classical space
of functions / distributions containing the solutions to \eqref{e:Phi4bis}.
The idea then is to enrich this space in a way that incorporates the minimal amount of information
required in order to be able to realise $\Phi \mapsto \Phi^3$ as a continuous map, but in a way
that remains consistent with the classical interpretation of \eqref{e:potential}, so that
one has hope to be able to obtain the type of universality results mentioned above.

A simple, but surprisingly effective, idea introduced in this context in \cite{DPD1,DPD2} (but see
also \cite{Bourgain} for the same idea in a slightly different context) is to try to solve \eqref{e:Phi4bis}
by setting $\Phi = \Psi + X$, where $X$ is given by \eqref{e:SHEadd}, so that it remains to solve
\begin{equ}[e:DPD]
\d_t \Psi = \Delta \Psi + c (\Psi + X) - (\Psi^3 + 3X \Psi^2 + 3X^2 \Psi + X^3)\;.
\end{equ}
The point here is that $\xi$ no longer appears explicitly in the right hand side,
so we can expect $\Psi$ to be more regular than $\Phi$, hopefully function-valued
so that the terms $\Psi^2$ and $\Psi^3$ are well-defined.

Assume now that $X$ is scaling invariant with some negative scaling exponent 
$-\alpha$ (in dimension $d$, this is the case for the solution to \eqref{e:SHEadd}
with $\alpha = {d\over 2} -1$). Then, even though we don't know a priori how to define
its powers, one would expect any reasonable definition of $X^p$ to be a stationary process
which is scaling invariant with exponent $-\alpha p$. 
Under relatively weak moment assumptions, such processes belong to the negative H\"older spaces
$\CC^\beta$ for all $\beta < -\alpha p$, which suggests that the least regular term on
the right hand side of \eqref{e:DPD} is the term $X^3$, with regularity just below $-3\alpha$,
so one would expect by Schauder's estimate $\Psi$ to have regularity just below $2-3\alpha$,
which is positive as soon as $\alpha < 2/3$. Since the product of a function in $\CC^\gamma$ with
a distribution in $\CC^\beta$ is well-defined (in the sense that there is a unique continuous 
bilinear extension of the usual product) if and only if $\beta + \gamma > 0$, this shows
that, \textit{assuming} that $X^2$ and $X^3$ are known, \eqref{e:DPD} is well-posed
as soon as $2-5\alpha > 0$, i.e.\ $\alpha < 2/5$.

In dimension $d=2$, the case treated in \cite{DPD2}, this is indeed the case. However,
it is still not clear how the processes $X^2$ and $X^3$ should be defined in that case. 
A ``naive'' construction of these processes would be to try to consider a regularisation
$X_\eps$ of $X$ and to check whether $X_\eps^2$ and $X_\eps^3$ have a limit. 
Write for example $\xi_\eps$ for a centred Gaussian random field with covariance given by 
$\rho_\eps(t,x) = \eps^{-3}\rho(t/\eps^2, x/\eps)$ for some compactly supported function $\rho$ integrating to $1$
and define $X_\eps$ by \eqref{e:SHEadd} with $\xi$ replaced by $\xi_\eps$. 
It can then be verified that 
\begin{equ}
\E X_\eps^2 = \int P(t,x)P(s,y)\,\rho_\eps(t-s,y-x)\,dt\,ds\,dx\,dy\;,
\end{equ}
where $P$ is the standard heat kernel. A simple calculation shows that,
in dimension $d=2$, one has $\E X_\eps^2 \sim \log \eps$, while in dimension $d=3$
one has $\E X_\eps^2 \sim \eps^{-1}$, so that there is no hope that $X_\eps^2$ converges to
a limiting proces as $\eps \to 0$.

However, setting $C_\eps = \E X_\eps^2$, it turns out that the Wick powers 
$X_\eps^{\diamond p} \eqdef H_p(X_\eps,C_\eps)$, where $H_p$ denotes the $p$th Hermite
polynomial do converge to a limit for $p$ small enough.
Indeed, writing $G_\eps$ for the covariance of $X_\eps$, namely
\begin{equ}
G_\eps(t,x) = \E X_\eps(0,0)X_\eps(t,x)\;,
\end{equ}
one can verify that the covariance of $X_\eps^{\diamond p}$ is simply given by 
$G_\eps^p$. If $X$ is self-similar with exponent $-\alpha$, then its covariance is homogeneous
of order $-2\alpha$, so that $G_\eps^p$ does indeed converge to $G^p$ in distribution,
provided that $2p\alpha < d+2$. (The reason why we have $d+2$ appearing here is that
this is precisely the threshold at which a homogeneous function with parabolic scaling
loses integrability at the origin.)
Recalling that $\alpha = {d\over 2}-1$, we conclude that $X_\eps^{\diamond 2}$
and $X_\eps^{\diamond 3}$ do indeed converge to non-trivial limits $X^{\diamond 2}$
and $X^{\diamond 3}$ in dimension $d < 4$.

In dimension $d=2$, it is therefore natural to define $\Phi_\eps = \Psi_\eps + X_\eps$
with $\Psi_\eps$ defined as
in \eqref{e:DPD}, but with $X^2$ and $X^3$ replaced by $X_\eps^{\diamond 2}$
and $X_\eps^{\diamond 3}$. Standard continuity estimates then show that $\Phi_\eps$ converges to
some limit $\Phi$, which we would like to interpret as a ``solution'' to \eqref{e:Phi4bis}.
This raises the question of the interpretation of the
resulting process, in particular whether it is possible to obtain it as a limit of
solutions to stochastic PDEs of the type	 \eqref{e:Phi4bis} with the noise replaced
by a smoothened version. 
For this, note that $\Psi_\eps$ solves
\begin{equs}
\d_t \Psi_\eps &= \Delta \Psi_\eps + c (\Psi_\eps + X_\eps) - (\Psi_\eps^3 + 3X_\eps \Psi_\eps^2 + 
3X_\eps^{\diamond 2} \Psi_\eps + X_\eps^{\diamond 3}) \label{e:actionSmooth}\\
&= \Delta \Psi_\eps + c (\Psi_\eps + X_\eps) - \big(\Psi_\eps^3 + 3X_\eps \Psi_\eps^2 + 
3(X_\eps^{2} - C_\eps) \Psi_\eps + (X_\eps^{3} - 3C_\eps X_\eps)\big)\\
&= \Delta \Psi_\eps + (c+3C_\eps) (\Psi_\eps + X_\eps) - (\Psi_\eps^3 + 3X_\eps \Psi_\eps^2 + 
3X_\eps^{2} \Psi_\eps + X_\eps^{3})\;.
\end{equs}
This shows that, for any fixed $\eps > 0$, $\Phi_\eps$ is indeed nothing but the solution to
\eqref{e:Phi4bis} with $\xi$ replaced by $\xi_\eps$ and $c$ replaced
by $c+3C_\eps$.

\subsection{Renormalisation}
\label{sec:renorm}

Let us step back for a moment and try to give a high-level overview of what's happening
here. We started from a formal expression \eqref{e:Phi4bis} which hints at the existence of
a family of processes $\Phi^c$ indexed by a parameter $c$. We then consider 
some natural approximation $\Phi_\eps^c$ of these processes and study it for small $\eps$. 
Unfortunately, for every fixed $c$, the sequence $\Phi_\eps^c$ diverges (or rather converges 
to something trivial in this particular case \cite{Ryser})
as $\eps \to 0$. However, the \textit{family} of processes 
$\{\Phi_\eps^c\}_{c \in \R}$ \textit{does} converge to a limiting family 
$\{\Phi^c\}_{c \in \R}$ as $\eps \to 0$. This is so because, when considering the whole 
family of processes as our object of interest, we do not care about the way in which it
is parametrised. We can therefore just as well reparametrise it by viewing it 
as the family $\{\Phi_\eps^{c+C_\eps}\}_{c \in \R}$
with $C_\eps$ as in the previous section. 

This idea first arose in theoretical physics in the earlier part of the 20th century
during the development of quantum field theory. The situation there is similar:
one starts from a family of formal expressions for the Lagrangian of a system and then
tries to use this in order to compute scattering amplitudes. 
These amplitudes are expressed as an infinite series of terms, each of which is 
described by a Feynman diagram. Some of these diagrams happen to encode diverging
integrals and are therefore not well-defined. The procedure that was eventually devised
in order to ``cure'' these divergencies can be summarised as follows:
\begin{itemize}
\item Construct a family $M_\eps^c$ of theories indexed by the same constants $c$ 
appearing in the original formal Lagrangian, as well as on some cut-off parameter $\eps$.
These theories are such that they give rise to well-defined Feynman diagrams for any fixed
$\eps > 0$ and they formally appear to approximate to the desired theory as $\eps \to 0$.
\item Find an $\eps$-dependent reparametrisation $g^\eps$ of the parameter space 
in such a way that, for every $c$, $M_\eps^{g^\eps c}$ converges to a non-trivial limiting 
theory as $\eps \to 0$.
\item Show that, modulo reparametrisation of its parameter space, the family of 
limits obtained in this way does not depend on the regularisation procedure 
used in the first step.
\end{itemize}

Whenever we are in such a situation, we call the theories $\hat M^c$ and $\hat M_\eps^c$
given by
\begin{equ}
\hat M^c = \lim_{\eps \to 0}\hat M_\eps^{c}\;,\qquad \hat M_\eps^{c} = M_\eps^{g^\eps c}\;,
\end{equ}
the ``renormalised theories''. In the context of QFT, the arguments $c$ parametrising the renormalised
theories are called the ``effective coupling constants'', while the corresponding arguments
$g^\eps c$ for the regularised theory are called the ``bare coupling constants''.

This turns out to be a rather generic state of affairs which also applies
\textit{mutatis mutandis} in the context of semilinear stochastic PDEs, where 
the $M^c$'s are solution theories to families of such equations.

\subsection{A simple example}

Before we proceed, let us now give a very simple and self-contained example of renormalisation. 
In this example, which 
does not at all have the pretense to be connected to physical reality in any sense,
a ``theory'' $M$ is given by a Schwartz distribution, i.e.\ a continuous linear map $M \colon \CS \to \R$,
where $\CS$ denotes the space of Schwartz functions.
Assume now that we would like to describe the family of Schwartz distributions formally given by 
$M^c(x) = {c_1\over |x|} - c_2 \delta_0(x)$. This is of course 
nonsensical since, for every test function $\phi$ with $\phi(0) \neq 0$, one would
then have
\begin{equ}[e:nonsense]
M^{c}(\phi) = \int_\R c_1 {\phi(x)\over |x|}\,dx - c_2 \phi(0) = \infty\;,
\end{equ}
so that $M^{c}$ really isn't a Schwartz distribution at all!

Given a continuous function $\eta \colon [-1,1] \to \R$ with $\eta(-1) = \eta(1) = 1$,
we have a natural way of regularising this theory by replacing the singular function $1/|x|$ by
the continuous function given by
\begin{equ}
{1 \over |x|_\eps} := 
\left\{\begin{array}{cl}
	1/|x| & \text{if $|x| \ge \eps$,} \\
	\eps^{-1}\eta(x/\eps) & \text{otherwise.}
\end{array}\right.
\end{equ}
This regularisation procedure then yields a family $M^{c,\eta}_\eps$ 
of well-defined  Schwartz distributions by simply replacing $|x|$ by
$|x|_\eps$ in \eref{e:nonsense}. In order to obtain finite quantities in the limit
$\eps \to 0$, we then set
\begin{equ}
c_1 = \hat c_1\;,\qquad
c_2 = \hat c_1 \int_{-1}^1 {dx \over |x|_\eps} + \hat c_2 = \hat c_1(\text{$\int\!\eta(x)dx$} - 2\log \eps) + \hat c_2\;.
\end{equ}
In this way, we have
\begin{equ}
M^{c(\eps,\eta,\hat c),\eta}_\eps(\phi) = 
\hat c_1  \int_{|x| \ge 1} {\phi(x)\over |x|}\,dx + \hat c_1 \int_{|x| < 1} {\phi(x) - \phi(0)\over |x|_\eps}\,dx - \hat c_2 \phi(0)\;.
\end{equ}
It is now immediate that this expression does indeed have a limit as $\eps \to 0$ and that
this limit is indeed independent of the choice of regularisation function $\eta$:
\begin{equ}
\hat M^{\hat c}(\phi) = 
\hat c_1  \int_{\R} {\phi(x) - \one_{|x| < 1} \phi(0)\over |x|}\,dx - \hat c_2 \phi(0)\;.
\end{equ}
(But for this we had to adjust the ``bare'' parameter $c_2$ in an $\eta$-dependent way!)

\subsection{Renormalisation of SPDEs}
\label{sec:renormSPDE}

Let us now give an overview of how the procedure explained in Section~\ref{sec:renorm} below
is implemented in the context of stochastic PDEs.
Consider a family of equations of the type
\begin{equ}[e:genSPDE]
\d_t u_i = L u_i + F_i^{(0)}(u) + \sum_{j=1}^m F_i^{(j)}(u)\,\xi_j\;,
\end{equ}
where $L$ is a strictly elliptic differential operator with constant coefficients
(for example $L = \Delta$ or $L = -\Delta^2$), the $\xi_j$ are finitely many stationary random distributions
which, at small scales, are close to scale invariant with exponents $\alpha_j$,
and the $F_i^{(j)}$ are local functions of $u$, in the sense that $F_i^{(j)}(u)(t,x)$ only
depends on $u$ and its derivatives evaluated at $(t,x)$.
We furthermore assume that this system of equations is \textit{locally subcritical}.
We do not give a precise definition of what this means, loosely speaking it means that 
the scaling behaviour of the non-linear terms appearing on the right hand side
of \eqref{e:genSPDE} is dominated at small scales by that of the noise terms, in the sense of
formal powercounting. It is a purely combinatorial condition, which only depends on the 
order of $L$, the exponents $\alpha_j$, and the degrees of the nonlinearities $F_i^{(j)}$
(with the convention that arbitrary smooth functions count as polynomials of infinite degree).

For example, in the case of the \eqref{e:Phi4}, the noise term in the right hand side
is $\xi$, which is scale invariant (under parabolic scaling) with exponent $-{d+2\over 2}$ in
the sense that
\begin{equ}[e:wnscale]
\eps^{\alpha}\xi(t/\eps^2, x/\eps) \eqlaw \xi(t,x)\;,\qquad \alpha = -{d+2\over 2}\;.
\end{equ}
On the other hand, the
``worst'' nonlinear term is $\Phi^3$ which is expected to be dominated at small scales
by a scale-invariant behaviour with exponent $-3\beta$, where $\beta = {d-2\over 2}$.
This is because the solution to $\d_t \Phi = \Delta \Phi + \xi$ displays is self-similar of 
order $\beta$ as a consequence of \eqref{e:wnscale} and self-similarity exponents do add up when
multiplying terms. (Or at least they would if these were functions that we are allowed to just multiply
together!)

We see that $-3\alpha > -{d+2\over 2}$ if and only if $d < 4$, so that this equation is
locally subcritical in dimensions below $4$, and we say that its ``critical dimension'' is $d=4$.
Note that this terminology is to some extent compatible with the one found in the physics literature, according
to which the critical dimension of the Ising model is also $4$. Performing the analogous
powercounting for the KPZ equation shows that it is locally subcritical for $d < 2$.
The following is a synthesis of the main results of \cite{Regularity,BHZ,Ajay,BCCH}:

\begin{theorem}\label{theo:main}
Consider \eqref{e:genSPDE} on a bounded $d$-dimensional torus and, for any 
smooth noise $\xi^{(\eps)}$, write $\CM(F,\xi^{(\eps)})$ for the classical solution 
to \eqref{e:genSPDE} with locally subcritical nonlinearity $F$ and noise $\xi$.

Then, there exists a finite-dimensional Lie group $\RR$ acting on the space of nonlinearities
as well as a map $\xi \mapsto g^{\xi} \in \RR$ defined on centred smooth stationary noises 
satisfying suitable moment conditions such that the map
\begin{equ}
(F, \xi) \mapsto \CM(g^\xi F,\xi)\;,
\end{equ}
extends continuously to all noises compatible with the scaling exponents $\alpha_j$.
\end{theorem}

The formulation of this result is somewhat imprecise on purpose. A precise formulation 
requires giving a formal definition of local subcriticality, a precise definition
of a class of admissible noises, as well as a topology in which they can be approximated.
This is surprisingly subtle and technical, so we avoid giving these details here.
For the same reason, we also ignored the question of a suitable class of initial conditions.
In the cases discussed in these notes, one can find $\alpha \in \R$ such that 
solutions exists for all $u_0 \in \CC^\alpha$, and such that these solutions then also take
values in $\CC^\alpha$ for positive times. This is however not the case in general: there 
are situations in which solutions
can only be defined for initial conditions that are ``sufficiently close to equilibrium'' in
a sense that can be made precise. This is the case for example when considering a model
analogous to \eqref{e:Phi4}, but in ``effective dimension'' $d = 4-\delta$ for small enough.
One way to set up such a model is to consider \eqref{e:Phi4} in dimension $4$, but to replace $\xi$ by
space-time white noise convolved with a kernel of the type $(t,x) \mapsto (t^2 + |x|^4)^{\delta-6 \over 4}$,
where the number $6$ refers to the scaling dimension of $4+1$-dimensional parabolic space-time. 
Finally, this general result is local in time, i.e.\ we work in a topology that allows solutions 
to blow up in finite time and all approximation results are only relevant up to this blow-up time.

\subsection{Structure of proof}
\label{sec:proof}

Let us now give a short overview of how a result like Theorem~\ref{theo:main}
can be proven. The first step, which is based on an idea first introduced by Lyons in
the nineties in the context of controlled ODEs (and in particular stochastic 
differential equations, see \cite{Terry,TerryBook1,TerryBook2,Peter1,Peter2}) is to 
factorise the classical solution map 
$\CM$ into a map that acts only on the noise and ``enhances'' it with additional 
information and a map that makes uses of this additional information to build the solution.
In the context of ODEs, this enhancement of the noise is called a ``rough path'', while in
the present context it is called a ``model'' for reasons that will become clearer later on.

For definiteness, let us write $\CN$ for a suitable space of admissible nonlinearities
leading to locally subcritical problems
and $\CY$ for a suitable space of distributions including possible blow-up times
arising in Theorem~\ref{theo:main}, so that the classical solution map
$\CM$ can be viewed as a map $\CM\colon \CN \times \CB_\infty \to \CY$, where the second factor
$\CB_\infty$ is a space of smooth noises.
We also write $\CB$ for the space of all distributions $(\xi_i)_i$ such that
$\xi_i \in \CC^{\alpha_i-\kappa}$ for some (small) $\kappa > 0$ and every $i \in \{1,\ldots,m\}$, 
with $\alpha_i$ as in \eqref{e:genSPDE}. We also view $\CB_\infty$ as a subspace of $\CB$.
From now on we completely ignore the additional dependence of the solution on the initial condition
since this does not introduce any additional conceptual difficulty.

The first step mentioned above is thus to find a topological space $\CX$ and
maps $\CL\colon \CB_\infty \to \CX$,  $\pi\colon \CX \to \CB$ and $\CM_A\colon \CN \times \CX \to \CY$
such that one has the identities
\begin{equ}[e:wanted]
\CM(F,\xi) = \CM_A(F, \CL \xi)\;,\qquad \pi \circ \CL = \id\;.
\end{equ}
Without additional constraints, it is of course trivial to do so: just take $\CX = \CB_\infty$
and $\CL$ and $\pi$ to be the identity. It is however highly non-trivial  to do this if we
impose the following additional constraints:
\begin{enumerate}
\item The map $\CM_A$ is \textit{continuous}.
\item The map $\pi$ is continuous and \textit{surjective}.
\end{enumerate}
Indeed, these conditions are competing: the smaller the space $\CX$, the easier it is
for $\CM_A$ to be continuous, but it needs to be large enough for $\pi$ to be surjective.
As a matter of fact, it was shown in \cite{TerryCounterexample} that already in very simple cases
it is impossible for any Banach space $\CX$ to simultaneously satisfy all these constraints!

It was shown in \cite{Regularity} (see also \cite{BHZ} for a cleaner formulation that allows for systems
of equations and collections of noises with different scaling exponents) that it is indeed 
possible to find such a space $\CX$
in all locally subcritical situations. The problem is that the map $\CL$, although continuous
on the space of smooth functions, does of course not extend in a continuous way to 
all of $\CB$, so that this is not yet sufficient to give a canonical interpretation to 
solutions to \eqref{e:genSPDE}: we would like to find a ``canonical'' random variable
$Z$ with values in $\CX$ such that the law of $\pi Z$ is equal to that of the noise $\xi$,
so that we could then \textit{define} solutions to \eqref{e:genSPDE} to be given by
$\CM_A(F,Z)$. Ideally, one would even want these random variables to be defined on the same 
probability space, so that $Z = \hat \CL(\xi)$ for some measurable map $\hat \CL$
satisfying $\pi \circ \hat \CL = \id$, which would then yield a corresponding
notion of ``strong solution''.

It was shown in \cite{Regularity} that there exists a finite-dimensional nilpotent 
Lie group $\RR$ acting on $\CX$ in a way that leaves $\pi$ invariant in the sense that
$\pi(gZ) = \pi(Z)$ for every $Z \in \CX$ and every $g \in \RR$. It is therefore natural
to try to build the random variable $Z$ by considering a regularisation $\xi^{(\eps)}$
of the noise $\xi$ and setting $Z^{(\eps)} = g^\eps \xi^{(\eps)}$. The
question is then whether there exist elements $g^\eps$ such that the $Z^{(\eps)}$ converge.
In \cite{Regularity}, this was shown for two examples on a case-by-case basis, but a
general theory was lacking due to a lack of understanding of the group $\RR$ and
its action on $\CX$.

This was partly remedied in \cite{BHZ}, where an explicit description of $\RR$ (or rather
a ``sufficiently large'' subgroup thereof) is given.
This allows \cite{Ajay} to show the following:

\begin{proposition}
In the context of Theorem~\ref{theo:main},
there exists a map $\xi \mapsto g^{\xi} \in \RR$ defined on centred smooth stationary noises 
satisfying suitable moment conditions such that the map
\begin{equ}
\xi \mapsto g^\xi \CL(\xi)\;,
\end{equ}
extends continuously to all noises compatible with the scaling exponents $\alpha_j$.
\end{proposition}

\begin{remark}
Note that this is a \textit{probabilistic} statement: it is important to consider random 
noises $\xi$ that are stationary in space-time. The continuity mentioned in the statement is not
a continuity at the level of samples of this process, but at the level of its law.
When considering convergences of the type $\rho_\eps \star \xi \to \xi$ for some 
mollifier $\rho_\eps$, the whole sequence naturally lives on the same probability space and the
result of \cite{Ajay} actually yields convergence in probability.
\end{remark}

\begin{remark}\label{rem:Xhat}
The map $\xi \mapsto g^\xi$ is constructed as follows. Given any class of locally subcritical
stochastic PDEs as described above, one first builds a linear space $\hat \CX$ as well as a 
polynomial (in the sense of ``finite sum of multilinear'') 
map $\hat \CL \colon \CB_\infty \to \hat \CX$ which is such that it maps stationary processes to stationary
processes under a suitable action of the group of translations on $\hat \CX$. The actual (nonlinear) space $\CX$ 
can then be viewed as a subset of $\hat \CX$ which turns into a Polish space when
endowed with a suitable metric. This metric is \textit{not} derived from any norm on $\hat \CX$ but instead encodes
the nonlinear features of $\CX$. The map $\CL$ mentioned in the proposition then coincides with $\hat \CL$
and the action of $\RR$ on $\CX$ is actually derived from an action on $\hat \CX$.
The element $g^\xi$ is then the unique element of $\RR$ with the property that 
\begin{equ}[e:center]
\E \bigl(g^\xi \hat \CL(\xi)\bigr)(0) = 0\;.
\end{equ}
The fact that such an element exists, is unique, that the action of $\RR$ on $\hat \CX$ leaves $\CX$
invariant and that it is continuous in the topology of $\CX$ are all non-trivial statements
that are contained in \cite{Regularity} for special cases and \cite{BHZ} in full generality.
\end{remark}

When combining this with the continuity of the map $\CM_A$, this shows that 
the map $\hat \CM$ defined by
\begin{equ}[e:dualAction]
\hat \CM(F,\xi) = \CM_A\bigl(F, g^\xi \CL(\xi)\bigr)\;,
\end{equ}
does indeed have the type of continuity properties with respect to the noise that are
announced in Theorem~\ref{theo:main}. The last missing ingredient in the proof
is the fact that $\hat \CM(\cdot,\xi)$ is really nothing but a reparametrisation of
$\CM(\cdot,\xi)$. This is shown in full generality in \cite{BCCH}:

\begin{proposition}
In the context of Theorem~\ref{theo:main},
there exists an action of $\RR$ onto $\CX$ on the right such that, for every $Z \in \CX$ and
every $F \in \CN$, 
the identity
\begin{equ}[e:renEqu]
\CM_A\bigl(F, g Z\bigr) =  \CM_A\bigl(Fg, Z\bigr)\;.
\end{equ}
holds for every $g \in \RR$. In particular, it follows that the renormalised solution map $\hat \CM$
is related to the classical solution map by $\hat \CM(F,\xi) = \CM(Fg^\xi,\xi)$. 
\end{proposition}

\section{Regularity structures}

In this final section, we give a description of the main ingredients appearing in the sketch of proof
given in the previous section, namely the space $\CX$ and the group $\RR$. 

\subsection{Revisiting the \texorpdfstring{$\Phi^4$}{Phi4} equation}
\label{sec:Phi4bis}

Let us first note that the construction given in Section~\ref{sec:Phi4} is precisely 
of the type described in the previous section with the important feature that in this case the space $\CX$
itself is linear. Indeed, we can set 
\begin{equ}
\CX = \CC^{-\alpha} \oplus \CC^{-2\alpha} \oplus \CC^{-3\alpha}\;,
\end{equ}
as well as 
\begin{equ}
\CL(\xi) = \bigl(X, X^2, X^3\bigr)\;,
\end{equ}
where $X$ is given as the solution to $\d_t X = \Delta X + \xi$. (We ignore the minor subtleties arising 
from the fact that such a solution may not necessarily exist. One way of circumventing this problem
is to actually set $\d_t X = (\Delta-1) X + \xi$ which leads to minor modifications in the formulas given below.)
Conversely, we set
\begin{equ}
\pi(X_1,X_2,X_3) = \d_t X_1 - \Delta X_1\;.
\end{equ}
The ``solution map'' $\CM_A$ is then given by 
postulating that, for $Z = (X_1,X_2,X_3)$,  $\Phi = \CM_A(c, Z)$ is given by 
$\Phi = \Psi + X_1$, with $\Psi$ solving
\begin{equ}[e:DPDbis]
\d_t \Psi = \Delta \Psi + c (\Psi + X_1) - (\Psi^3 + 3X_1 \Psi^2 + 3X_2 \Psi + X_3)\;.
\end{equ}
Comparing this to \eqref{e:DPD}, we immediately conclude that \eqref{e:wanted} does 
indeed hold.

Let now $\RR = (\R^2, +)$ and, for $g \in \R^2$, define its action on $\CX$ by
\begin{equ}[e:actionPhi42]
gZ = (X_1, X_2 - g, X_3 - 3g X_1)\;.
\end{equ}
The calculation in \eqref{e:actionSmooth} then shows that the identity \eqref{e:dualAction}
does hold, provided that the action of $\RR$ onto the space of right hand sides
(which in this case is parametrised by the constant $c$) is given by
$g \colon c \mapsto c + 3g$.
The map $\xi \mapsto g^\xi$ is then given by
\begin{equ}
g^\xi = \E X^2(0)\;,\qquad \d_t X = \Delta X + \xi\;.
\end{equ}
If we restrict ourselves to symmetric noises, or even just to noises with vanishing first and third moments,
it immediately follows that this choice of $g^\xi$ does indeed satisfy \eqref{e:center}.
The fact that it also makes $\xi \mapsto g^\xi \CL(\xi)$ continuous
is less obvious and relies on the particular form of the action \eqref{e:actionPhi42}:
the constant $3$ appearing in the last component seems somewhat arbitrary, but it serves to 
seemingly unrelated purposes. First, if we restrict ourselves to actions of the form
\eqref{e:actionPhi42}, the value $3$ for this constant is the only possible choice
which guarantees that there exists a dual action of $\RR$ on the space of right hand sides for
our equation such that \eqref{e:dualAction} holds. The second consequence of this particular choice
is that $\xi \mapsto g^\xi \CL(\xi)$ is continuous provided that we restrict ourselves to
stationary processes satisfying a natural cumulant bound.

Although we will not give a proof of this latter fact here, let's perform a calculation 
that makes it plausible, and that simultaneously shows what kind of probabilistic assumptions
one should impose on the noise $\xi$. Fix a test function $\phi$ and consider the random variables
\begin{equ}
\hat Y_\phi = \int \phi(z) X(z)\,dz\;,\qquad 
Y_\phi = \int \phi(z) \bigl(X^2(z) - \E X^2(z)\bigr)\,dz\;,
\end{equ}
with $X$ given by \eqref{e:SHEadd}. 
Since $X$ is stationary, $\E X^2(z) = \E X^2(0)$ and these are precisely the random variables
obtained by testing the first two components of $g^\xi \CL(\xi)$ with $\phi$.

We now give an example in dimension $3$ of a sequence of 
non-Gaussian processes $\xi_\eps$ 
so that the corresponding random variables $Y_\phi^{(\eps)}$ and $\hat Y_\phi^{(\eps)}$ converge to a limit, but so that
 $\E X_\eps^2(z)$ diverges. For this, let $\eta$ be a stationary generalised Gaussian random
 field with covariance given by $\E \eta(z)\eta(z') = K(z-z') \sim |z-z'|^{-{9\over 4}}$, 
with $|z|$ denoting the parabolic norm of $z = (t,x)$. We then set $\eta_\eps = \rho_\eps \star \eta$,
write $K_\eps$ for the covariance function of $\eta_\eps$, and set
$\xi_\eps(z) = \eta_\eps^2(z) - K_\eps(0)$. 

We now introduce a graphical notation for integrals involving $K_\eps$, the heat kernel $P$, and
the test function $\phi$. We denote integration variables by nodes of our graphs, with the origin
drawn as a special node \tikz[baseline=-3] \node [root] {};. 
Each line then represents a kernel, with 
\tikz[baseline=-0.1cm] \draw[kernel] (0,0) to (1,0);
representing the heat kernel $P$, 
\tikz[baseline=-0.1cm] \draw[rho] (0,0) to (1,0);
representing the covariance $K_\eps$, and
\tikz[baseline=-0.1cm] \draw[testfcn] (1,0) to (0,0);
representing the test function $\phi$. The argument of a kernel is the difference between the
values of the variables represented by the corresponding target and base nodes.
The line representing $K_\eps$ are not oriented because this kernel is symmetric, being the covariance
of a stochastic process.
With this notation, setting for example
\begin{equ}
\hat Y_\phi = \int \phi(z) X(z)\,dz\;,
\end{equ}
it follows from the definition of $X$ and Wick's formula that 
\begin{equ}
\E \hat Y_\phi^2 = 2
\int \phi(z_1)\phi(z_2)P(z_1-z_3)P(z_2-z_4) K_\eps(z_3-z_4)^2\,dz
= 2\;
\begin{tikzpicture}[scale=0.35,baseline=0.5cm]
	\node at (0,0)  [root] (root) {};
	\node at (-1.5,1)  [dot] (left) {};
	\node at (-1.5,3)  [dot] (left1) {};
	\node at (1.5,1)  [dot] (right) {};
	\node at (1.5,3)  [dot] (right1) {};
	
	\draw[testfcn] (left) to (root);
	\draw[testfcn] (right) to (root);
	
	\draw[kernel] (left1) to (left);
	\draw[kernel] (right1) to (right);
	\draw[rho,bend left=60] (left1) to (right1); 
	\draw[rho,bend left=60] (right1) to (left1); 
\end{tikzpicture}\;.
\end{equ}
Note now that in dimension $3+1$ with parabolic scaling, if $K_i$ are kernels behaving like $K_i(z) \sim |z|^{-\alpha_i}$
(we say that they are self-similar of order $\alpha_i$)
and one has $\alpha_i < 5$ while $\alpha_1 + \alpha_2 > 5$, then 
$K_1 \star K_2$ is a kernel which is again self-similar of order $\alpha_1 + \alpha_2 - 5$.
In this particular example, $K_\eps^2$ is self-similar of order ${9\over 2}$, while the heat
kernel is self-similar of order $3$, so that the kernel
\begin{equ}
\bar P \star K_\eps^2 \star P = 
\begin{tikzpicture}[scale=0.35,baseline=-0.1cm]
	\node at (-1,0)  [dot] (left) {};
	\node at (1,0)  [dot] (right) {};
	
	\draw[kernel] (left) to (-3,0);
	\draw[kernel] (right) to (3,0);
	\draw[rho,bend left=60] (left) to (right); 
	\draw[rho,bend left=60] (right) to (left); 
\end{tikzpicture}\;,
\end{equ}
where $\bar P(z) = P(-z)$, is self-similar of order ${1\over 2}$ and therefore locally integrable.
A similar calculation allows to conclude that $\lim_{\eps \to 0} \E (\hat Y_\phi^{(\eps)} - \hat Y_\phi)^2 = 0$
as claimed.

For this, no renormalisation was necessary. The situation gets a bit more interesting regarding 
$Y_\phi$. It follows from Wick's formula combined with the definition of $\xi_\eps$ that one has
\begin{equ}[e:variance]
\E Y_\phi^2 = 8\;
\begin{tikzpicture}[scale=0.35,baseline=0.5cm]
	\node at (0,0)  [root] (root) {};
	\node at (-1.5,1)  [dot] (left) {};
	\node at (-1.5,2.5)  [dot] (left1) {};
	\node at (1.5,1)  [dot] (right) {};
	\node at (1.5,2.5)  [dot] (right1) {};
	\node at (-1.5,4.5)  [dot] (left2) {};
	\node at (1.5,4.5)  [dot] (right2) {};
	
	\draw[testfcn] (left) to (root);
	\draw[testfcn] (right) to (root);
	
	\draw[kernel] (left1) to (left);
	\draw[kernel] (right1) to (right);

	\draw[kernel,bend right=60] (left2) to (left);
	\draw[kernel,bend left=60] (right2) to (right);

	\draw[rho,bend left=60] (left1) to (right1); 
	\draw[rho,bend left=60] (right1) to (left1); 
	\draw[rho,bend left=60] (left2) to (right2); 
	\draw[rho,bend left=60] (right2) to (left2); 
\end{tikzpicture} +
32\;
\begin{tikzpicture}[scale=0.35,baseline=0.5cm]
	\node at (0,0)  [root] (root) {};
	\node at (-1.5,1)  [dot] (left) {};
	\node at (-1.5,2.5)  [dot] (left1) {};
	\node at (1.5,1)  [dot] (right) {};
	\node at (1.5,2.5)  [dot] (right1) {};
	\node at (-1.5,4.5)  [dot] (left2) {};
	\node at (1.5,4.5)  [dot] (right2) {};
	
	\draw[testfcn] (left) to (root);
	\draw[testfcn] (right) to (root);
	
	\draw[kernel] (left1) to (left);
	\draw[kernel] (right1) to (right);

	\draw[kernel,bend right=60] (left2) to (left);
	\draw[kernel,bend left=60] (right2) to (right);

	\draw[rho] (left1) to (right1); 
	\draw[rho] (right1) to (right2); 
	\draw[rho] (left2) to (right2); 
	\draw[rho] (left1) to (left2); 
\end{tikzpicture} +
16\;
\begin{tikzpicture}[scale=0.35,baseline=0.5cm]
	\node at (0,0)  [root] (root) {};
	\node at (-1.5,1)  [dot] (left) {};
	\node at (-1.5,2.5)  [dot] (left1) {};
	\node at (1.5,1)  [dot] (right) {};
	\node at (1.5,2.5)  [dot] (right1) {};
	\node at (-1.5,4.5)  [dot] (left2) {};
	\node at (1.5,4.5)  [dot] (right2) {};
	
	\draw[testfcn] (left) to (root);
	\draw[testfcn] (right) to (root);
	
	\draw[kernel] (left1) to (left);
	\draw[kernel] (right1) to (right);

	\draw[kernel,bend right=60] (left2) to (left);
	\draw[kernel,bend left=60] (right2) to (right);

	\draw[rho] (left1) to (right2); 
	\draw[rho] (right1) to (left1); 
	\draw[rho] (left2) to (right1); 
	\draw[rho] (right2) to (left2); 
\end{tikzpicture}\;.
\end{equ}
The main point here is that thanks to the fact that we subtracted $\E X^2(0)$ in the definition
of the second component, there is no term of the form 
\begin{equ}[e:bad]
\begin{tikzpicture}[scale=0.35,baseline=0.1cm]
	\node at (0,0)  [root] (root) {};
	\node at (-1.5,1)  [dot] (left) {};
	\node at (-3,1)  [dot] (left1) {};
	\node at (-5,1)  [dot] (left2) {};
	\node at (1.5,1)  [dot] (right) {};
	\node at (3,1)  [dot] (right1) {};
	\node at (5,1)  [dot] (right2) {};
	
	\draw[testfcn] (left) to (root);
	\draw[testfcn] (right) to (root);
	
	\draw[kernel] (left1) to (left);
	\draw[kernel] (right1) to (right);

	\draw[kernel,bend right=60] (left2) to (left);
	\draw[kernel,bend left=60] (right2) to (right);

	\draw[rho,bend right=80] (left1) to (left2); 
	\draw[rho,bend left=80] (right1) to (right2); 
	\draw[rho] (right1) to (right2); 
	\draw[rho] (left1) to (left2); 
\end{tikzpicture} 
\end{equ}
appearing in this expression. Instead, all the terms appearing in \eqref{e:variance}
are $2$-connected: one needs to cut at least two edges in order to disconnect them.
The study of the convergence of multiple integrals of this type is a standard ingredient
of perturbative quantum field theory. In this case, one can apply Weinstein's theorem
\cite{Weinberg} (see also \cite[Prop.~2.3]{BPHZ}) which can be formulated as follows.
Assign to each edge $e$ an exponent $a_e$ by setting $a_e= 0$ for edges representing $\phi$,
$a_e = 3$ for edges representing $P$, and $a_e = {9\over 4}$ for edges representing $K_\eps$.
Then, integrals described by a graph $\Gamma$ as above converge if, for every subgraph 
$(\bar V, \bar E) \subset \Gamma$ (with $\bar V$ a subset of the vertices of $\Gamma$ and $\bar E$
a subset of its edges such that the endpoints of every edge in $\bar E$ belongs to $\bar V$),
one has the bound
\begin{equ}[e:Weinstein]
\sum_{e \in \bar E} a_e < 5\bigl(|\bar V|-1\bigr) \;.
\end{equ}
One can verify that all the graphs appearing in \eqref{e:variance} do indeed have this property
while the graph \eqref{e:bad} does not. (Just take as a subgraph 
the graph \begin{tikzpicture}[scale=0.35,baseline=0.2cm]
	\node at (1.5,1)  [dot] (right) {};
	\node at (3,1)  [dot] (right1) {};
	\node at (5,1)  [dot] (right2) {};
	
	\draw[kernel] (right1) to (right);
	\draw[kernel,bend left=40] (right2) to (right);

	\draw[rho,bend left=40] (right1) to (right2); 
	\draw[rho] (right1) to (right2); 
\end{tikzpicture}
for which one has $\sum_{e \in \bar E} a_e = 10{1\over2} > 10$.)

What we've learned from this calculation is that the key property of the renormalisation
procedure \eqref{e:center} is that it creates cancellations for terms of lower
connectivity in the graphical representation for the variances of the random
variables $Y_\phi$ and $\hat Y_\phi$. In general, one has a decomposition of the
type \eqref{e:variance}, but with instances of
$2\begin{tikzpicture}[scale=0.35,baseline=-0.13cm]
	\node at (0,0)  [dot] (left) {};
	\node at (2,0)  [dot] (right) {};
	
	\draw[rho,bend left=40] (left) to (right);
	\draw[rho,bend left=40] (right) to (left); 
\end{tikzpicture}$
replaced by the covariance of $\xi_\eps$ and instances of 
$32\begin{tikzpicture}[scale=0.35,baseline=0.05cm]
	\node at (0,0)  [dot] (left) {};
	\node at (0,1)  [dot] (left1) {};
	\node at (2,0)  [dot] (right) {};
	\node at (2,1)  [dot] (right1) {};
	
	\draw[rho] (left) -- (right) -- (right1) -- (left1) -- (left);
\end{tikzpicture}
+
16\begin{tikzpicture}[scale=0.35,baseline=0.05cm]
	\node at (0,0)  [dot] (left) {};
	\node at (0,1)  [dot] (left1) {};
	\node at (2,0)  [dot] (right) {};
	\node at (2,1)  [dot] (right1) {};
	
	\draw[rho] (left) -- (right) -- (left1) -- (right1) -- (left);
\end{tikzpicture}$
replaced by its four-point cumulant. It is then natural to impose conditions
on these cumulants that are akin to what we observed above. The key condition
one imposes is that for $k = \ell + m$ and writing $\kappa_k$ for the
$k$-point cumulant, the behaviour of 
$\kappa_{\ell+m}(z_1,\ldots,z_\ell,z_1',\ldots,z_m')$ is less singular than 
that of $\kappa_{\ell}(z_1,\ldots,z_\ell)\kappa_{m}(z_1',\ldots,z_m')$
as some of the variables $z_i$ get close to each other. In other words, one
assumes that cumulants behave ``better'' than what one can deduce from their 
expressions in terms of moments. This is a very natural condition which, as we have
just seen in an example, is verified for very large classes of random fields.

\begin{remark}
One can also consider noises with non-vanishing third moments, 
but this then forces
us to consider a slightly larger class of right hand sides, for example
we can consider all equations of the type
\begin{equ}[e:nonSym]
\d_t \Phi = \Delta \Phi + c_0 + c_1 \Phi + c_2\Phi^2 - \Phi^3 + \xi\;.
\end{equ}
Furthermore, one should then take $\RR = (\R^2,+)$ acting on the space $\CX$ by
\begin{equ}
gZ = (X_1, X_2 - g_1, X_3 - 3g_1 X_1 - g_2)\;.
\end{equ}
In this case, it is a good exercise to show that 
the dual action on \eqref{e:nonSym} for which \eqref{e:dualAction} holds 
is given by
\begin{equ}
(c_0,c_1,c_2) \mapsto (c_0-c_2g_1+g_2,c_1 + 3g_1,c_2)\;.
\end{equ}
\end{remark}

\subsection{Basic definitions}

In the above example, the space $\CX$ is a Banach space and the group $\RR$ is
abelian. At the level of generality considered in Theorem~\ref{theo:main}, it is
typically not possible to enforce this and a more sophisticated approach is required.

The problem is that in general, the trick of subtracting some fixed process from 
the solution in order to improve its regularity properties does not work. 
The idea then is to perform a similar procedure, but to proceed \textit{locally} rather than
\textit{globally} and to determine these local terms in a self-consistent way. 
More precisely, the aim is to provide a local description of the solution by a kind of 
Taylor expansion and to then find a fixed point problem for the coefficients of this expansion.
The usual Taylor polynomials won't do of course since, as we have already seen above,
solutions may not even be function-valued. 

Instead, as in the previous subsection, the idea is to build a collection of noise-dependent
objects (in this case $X$, $X^2$ and $X^3$)
which are useful in order to describe both the solution to \eqref{e:genSPDE} and
its right hand side. In general, these objects can be arranged into a structure that closely mimics
that of the usual Taylor polynomials. 
Before we proceed, we describe this structure in a very general context that does not 
refer to stochastic PDEs at all.

The starting point for our construction is a vector space $T$ that should be though of as containing the
coefficients of our ``Taylor-like'' expansion at any point. It is natural to 
postulate that $T$ is a graded vector space $T = \bigoplus_{\alpha \in A} T_\alpha$,
for some discrete set $A$ of possible ``homogeneities''. 
For example, in the case of the usual Taylor expansions,
we take for $A$ the set of natural numbers and $T_\ell$ contains the coefficients corresponding
to the monomials of total degree $\ell$. In general, we only assume that the set $A$ is bounded from
below and locally finite, and that each $T_\alpha$ is a 
real Banach space, although in many examples of interest these spaces will be finite-dimensional.

A crucial characteristic of
Taylor expansions is that an expansion around any point $x_0$ can be 
re-expanded around any other point $x_1$, namely simply by making use of the identity
\begin{equ}[e:TaylorExp]
(x-x_0)^m = \sum_{k+\ell = m} \binom{m}{k} (x_1 - x_0)^k\cdot (x-x_1)^\ell\;.
\end{equ}
In the general case, we only assume that there are linear maps $\Gamma_{xy}$
transforming the coefficients of an expansion around $y$ into the coefficients for the 
same ``polynomial'', but this time expanded around $x$.

In view of the example of Taylor expansions, it is natural to impose that any such ``reexpansion map''
$\Gamma_{xy}$ has the property that
if $\tau \in T_\alpha$, then $\Gamma_{xy} \tau - \tau \in \bigoplus_{\beta < \alpha} T_\beta
=: T_{<\alpha}$. In other words, 
when reexpanding a homogeneous monomial around a different point, 
the leading order coefficient remains 
the same, but lower order monomials may appear, just as is the case in \eqref{e:TaylorExp}. 
Furthermore, one should be able to compose reexpansions,
since taking an expansion around $x$, reexpanding it around $y$ and then reexpanding the result
around a third point $z$ should be the same as reexpanding the first expansion around $z$.
In other words, it seems natural to impose the identity $\Gamma_{xy}\Gamma_{yz} = \Gamma_{xz}$.
These considerations can be summarised in the following definition of an 
algebraic structure which we call a \textit{regularity structure}:

\begin{definition}\label{def:reg}
Let  
$T = \bigoplus_{\alpha \in A} T_\alpha$ be a vector space graded by $A \subset \R$ (discrete bounded below) 
such that each $T_\alpha$ is a 
Banach space. Let furthermore $G$ be a group of continuous operators on $T$ such that, for every $\alpha \in A$,
every $\Gamma \in G$, one has 
\begin{equ}[e:triang]
\tau \in T_\alpha\qquad\Rightarrow\qquad \Gamma \tau - \tau \in T_{<\alpha}\;.
\end{equ}
The pair $\TT=(T,G)$ is called a \textit{regularity structure} with \textit{model space} $T$ and \textit{structure group} $G$.
\end{definition}

\begin{remark}
We say that an element $\tau \in T_\alpha$ is ``homogeneous of degree $\alpha$'' and write
$\deg \tau = \alpha$.
\end{remark}

Such a regularity structure is a purely algebraic construct, but its purpose is to be 
endowed with some analytic ``flesh'': for each point $x$, we consider a linear map
$\Pi_x \colon T \to \CS'$ (here $\CS$ is the space of smooth test functions with compact support) 
with the property that, for every homogeneous $\tau$ of 
degree $\alpha$, the distribution $\Pi_x \tau$ satisfies an analytic bound that is homogeneous of degree $\alpha$ 
around $x$:
\begin{equ}[e:boundModel]
\bigl|\bigl(\Pi_x \tau\bigr)(\phi_x^\lambda)\bigr| \lesssim \lambda^\alpha\;,\qquad \forall \lambda \in (0,1]\;,
\end{equ}
where $\phi_x^\lambda$ denotes a test function $\phi$ of size $1$, translated so that its 
support contains $x$, and then rescaled so that its integral doesn't change but the diameter of
its support is of order $\lambda$. For a precise definition, see \cite[Def.~2.17]{Regularity}.

The reexpansion property mentioned above then suggests that we should restrict our attentions
to maps $\Pi$ as above with the additional property that one can find elements $\Gamma_{xy} \in G$ 
such that the identity 
\begin{equ}[e:alg]
\Pi_x \Gamma_{xy} = \Pi_y\;,
\end{equ}
holds for any pair of points $(x,y)$. By scaling arguments and by analogy with \eqref{e:TaylorExp}
it is furthermore natural to impose that if $\tau \in T_\alpha$ and $\beta < \alpha$, then the 
component of $\Gamma_{xy}\tau$ in $T_\beta$ should be of size at most $\CO(|x-y|^{\alpha-\beta})$.
We call such a pair $(\Pi, \Gamma)$ a ``model'' for the underlying regularity
structure. The space $\CX$ of all models is endowed with a natural topology which turns it into
a complete separable metric space, but not a linear space.
In the analogy with usual Taylor polynomials, the regularity structure encodes the algebraic properties
of Taylor monomials, while a model realises them as actual functions defined on some underlying Euclidean space.

\begin{example}\label{ex:poly}
The usual Taylor polynomials are cast in this framework as follows. Take for $T$ the space of
all polynomials $\sum_k a_k \sX^k$ in some abstract indeterminate $X$ (say there are $d$ indeterminates
so we interpret $k$ as a multiindex), with $T_m$ the subspace
spanned by those $\sX^k$ with $|k| = m$. The group $G$ is then isomorphic to $\R^d$,
acting on $T$ by $G_h \sX^k = (\sX-h)^k$ for $h \in \R^d$. The canonical model for this polynomial regularity
structure is given by setting
\begin{equ}
\bigl(\Pi_x \sX^k\bigr)(z) = (z-x)^k\;,\qquad 
\Gamma_{xy} = G_{y-x}\;.
\end{equ}
It is then easy to check that \eqref{e:boundModel} and \eqref{e:alg} do indeed hold.
\end{example}

We can then consider functions $f \colon \R^d \to T$ which are interpreted as
providing a ``jet'' $\Pi_x f(x)$ around every point $x$. Given a model, there 
are natural H\"older-type topology on the space of such functions which allow
one to define analogues $\CD^\gamma$ to the usual H\"older spaces $\CC^\gamma$.
The way these spaces are defined is completely analogous to the usual H\"older spaces
in the sense that one would like $f$ to be ``approximated by a polynomial up to order $\gamma$''.
In our context, a ``polynomial'' is naturally described by a function $p\colon \R^d \to T$
such that 
\begin{equ}[e:poly]
p(x) = \Gamma_{xy} p(y)\;,
\end{equ}
for any two points $x$ and $y$. This function can be
``reconstructed'' to a bona fide distribution on $\R^d$ by setting $\CR p = \Pi_x p(x)$.
Note that, by \eqref{e:alg} and \eqref{e:poly}, this is independent of the particular choice of $x$.
(This is of course not to say that the distribution $\CR p$ is constant!)

Given this discussion, it is natural to say that a function $f \colon \R^d \to T$ is in $\CD^\gamma$
if \eqref{e:poly} holds ``up to order $\gamma$''. More precisely, we impose that
\begin{equ}[e:Dgamma]
\|f(x) - \Gamma_{xy} f(y)\|_\alpha\lesssim |x-y|^{\gamma-\alpha}\;,\qquad \alpha < \gamma\;,
\end{equ}
where $\|\cdot\|_\alpha$ denotes the norm of the component in $T_\alpha$ of some element of $T$.
The reason for the exponent $\gamma - \alpha$ is that the component of $f$ in $T_\alpha$ should be
thought of as representing a type of ``derivative of order $\alpha$'', again by analogy with usual Taylor
expansions, so that \eqref{e:Dgamma} is analogous to the fact that if a function is of class $\CC^\gamma$,
then its derivative of order $\alpha$ is of class $\CC^{\gamma-\alpha}$, provided of course that $\alpha < \gamma$.
One fundamental result of the theory of regularity structures is the following, which 
is the analogue in this context of the ``sewing lemma'' of \cite{Max} and states that the operation $\CR$
defined above on polynomials extends canonically to all of $\CD^\gamma$ provided that $\gamma > 0$.

\begin{theorem}\label{theo:reconstr}
Given a regularity structure endowed with a model $Z = (\Pi,\Gamma)$, 
for any $f \in \CD^\gamma$ with $\gamma > 0$, there exists a unique distribution
$\CR f$ such that, for every $x \in \R^d$, the distribution $\Pi_x f(x) - \CR f$
vanishes at order $\gamma$ in the vicinity of $x$.
\end{theorem}

\begin{remark}
We will sometimes write $\CR(f,Z)$ instead to $\CR f$ to emphasise the fact that this
depends not only on the element $f \in \CD^\gamma$, but also on the underlying
 model $Z \in \CX$.
\end{remark}

\begin{remark}\label{rem:pointwise}
In the case when $\Pi_x\tau$ is actually a continuous function for every $\tau \in T$ one
has the explicit formula $\bigl(\CR f\bigr)(x) = \bigl(\Pi_x f(x)\bigr)(x)$. In general, 
the right hand side of this expression makes of course no sense since $\Pi_x f(x)$ might be 
a distribution which cannot be evaluated pointwise.
\end{remark}

Furthermore, the map $(Z,f) \mapsto \CR f$ mapping a ``model'' $Z = (\Pi,\Gamma)$ as well as
a ``modelled distribution'' $f \in \CD^\gamma$ onto the Schwartz distribution $\CR f$ is 
jointly continuous. (The precise continuity statement is somewhat subtle since the definition of the space $\CD^\gamma$
is itself dependent of the model $Z$, so that the pairs $(Z,f)$ really take values in a type of 
topological vector bundle.)
This theorem also justifies the terminology ``model'': the bound given in this theorem
states that the distribution $\CR f$ is modelled locally by $\Pi$ in the same way that a smooth
function is modelled locally by polynomials.

The link to the discussion in Section~\ref{sec:proof} is now the following.
The space $\CX$ is taken to be a suitable closed subset of the space of all models
for a regularity structure that is canonically associated with a given class of
semilinear parabolic stochastic PDEs. The map $\CM_A$ is then built by mimicking the
usual proof of well-posedness for parabolic PDEs, but in a suitable weighted version
of the space $\CD^\gamma$. This requires to build a whole calculus in these spaces in order to
give a meaning to the various operations appearing there.

\subsection{Calculus for regularity structures}
\label{sec:calculus}

Let us now show how one can construct the map $\CM_A \colon \CN \times \CX \to \CY$ alluded
to earlier. Recall that $\CN$ is a suitable space of right hand sides for the equation
of interest and $\CX$ will be chosen as the space of models for a suitable regularity structure.
The precise construction of the regularity structure in question is of course part of
the question. The idea is to formulate an SPDE of the type \eqref{e:genSPDE} as a fixed point problem in
some space $\CD^\gamma$ and to then set
\begin{equ}[e:defMA]
\CM_A(F,Z) = \CR(\hat \CM(F,Z),Z)\;,
\end{equ}
where $\hat \CM \colon \CN \times \CX \to \CD^\gamma$ is the solution to the fixed point problem
and $\CR$ is the reconstruction operator given by Theorem~\ref{theo:reconstr}.

In order to construct a regularity structure adapted to an equation of the type 
\eqref{e:genSPDE}, the idea is to start from the polynomial regularity structure endowed
with its canonical model as described in Example~\ref{ex:poly} and to systematically enlarge it
until it is sufficiently rich to support the operations required to formulate 
\eqref{e:genSPDE}. 
Our first step is to add to the indeterminates $\sX_i$ additional symbols $\sXi_i$
on which the structure group acts trivially. The model $(\Pi,\Gamma) = \CL(\xi)$ is then 
defined to be such that $\Pi_x \sXi_i = \xi_i$ for every $x$ and every $i$.
The degree $\alpha_i = \deg\sXi_i$ of $\sXi_i$ is constrained by the regularity of the generalised random fields
$\xi_i$ for which we want to build a solution theory in the sense that $\alpha_i$ should
be sufficiently small so that the bound \eqref{e:boundModel} holds for $\tau = \sXi_i$. In the case
when $\xi_i$ is white noise for example, this imposes the constraint $\alpha_i < -{D\over 2}$, where
$D$ is the scaling dimension of the underlying space(-time). 

It is natural to also add symbols of the type $\sX^k \sXi_i$ to the collection of basis vectors of 
$T$ and to impose that $\deg (\sX^k \sXi_i) = \deg \sX^k + \deg \sXi_i = \alpha_i + |k|$.
The natural action of $\R^d$ on these symbols is given by 
$G_h (\sX^k \sXi_i) = (\sX-h)^k \sXi_i$, and we naturally extend $\CL(\xi)$ to it
by setting $\bigl(\Pi_x \sX^k \sXi_i\bigr)(z) = (z-x)^k\xi_i(z)$.
It is an easy exercise to verify that these definitions are compatible with \eqref{e:boundModel} 
and \eqref{e:alg}.

\begin{remark}\label{rem:notalg}
In general, there is no reason to introduce symbols of the type $\sXi_i\sXi_j$
since products of noises do not naturally appear in \eqref{e:genSPDE}. Furthermore, in many situations
there is no natural way of renormalising the pointwise product $\xi_i \xi_j$ in a meaningful way.
As a consequence, our structure space $T$ will in general \textit{not} be an algebra. Although our notation
suggests that it is endowed with a product, this has a non-trivial domain of definition so that not
any two vectors can be multiplied with each other.
\end{remark}

At this stage, we note that we have enough structure to be able to use our theory to multiply
the distributions $\xi_i$ with sufficiently regular functions $f \in \CC^\gamma$. Indeed,
given $f \in \CC^\gamma$, we can canonically lift it to $F \in \CD^\gamma$ by setting
\begin{equ}[e:Taylor]
F(x) = \sum_{|k| < \gamma} {f^{(k)}(x) \over k!} \sX^k\;.
\end{equ}
Since $F(x)$ lives in the polynomial ``sector'' of our regularity structure, we can multiply
it by $\sXi_i$ and define $\hat F_i(x) = F(x)\sXi_i$.  
The condition $F \in \CD^\gamma$ implies that $\hat F_i \in \CD^{\gamma + \alpha_i}$ since
\begin{equ}
\hat F_i(x) - \Gamma_{xy}\hat F_i(y) = \bigl(F_i(x) - \Gamma_{xy}F_i(y)\bigr)\sXi_i\;,
\end{equ}
and multiplication by $\sXi_i$ changes the degree by $\alpha_i$, so that $\CR \hat F_i$ 
is well-defined as soon as $\gamma$ is sufficiently large so that $\gamma + \alpha_i > 0$.
This is nothing but the well-known fact that the product $(f,\xi) \mapsto f\cdot \xi$ extends
continuously to $\CC^\gamma \times \CC^\alpha$ if and only if $\gamma + \alpha > 0$, see
\cite{BookChemin}, so it seems that we haven't really learned anything new yet.

However, one verifies that this construction generalises very easily. Indeed, consider 
any two subspaces $V_1, V_2 \subset T$ that are invariant under the action of $G$, block-diagonal 
with respect to the decomposition $T = \bigoplus_\alpha T_\alpha$, and such that
there exists a bilinear ``product'' $\ast \colon V_1 \times V_2 \to T$ with the additional
property that, for every $\Gamma \in G$ and homogeneous $\tau_i \in V_i$, their product is homogeneous
and one has
\begin{equ}[e:algProd]
\Gamma(\tau_1\ast \tau_2) = (\Gamma \tau_1) \ast (\Gamma \tau_2)\;,\qquad
\deg (\tau_1\ast \tau_2) = \deg \tau_1 + \deg \tau_2\;.
\end{equ}
We furthermore set $\beta_i \le 0$ so that $V_i \subset \bigoplus_{\alpha \ge \beta_i} T_\alpha$.
One then has the following general result, the proof of which is an elementary exercise.

\begin{proposition}\label{prop:prod}
In the above context, given $F_i \in \CD^{\gamma_i}$ with values in $V_i$ and $\gamma_i > 0$,
one has $F_1 \ast F_2 \in \CD^\gamma$ with 
$\gamma = (\beta_1 + \gamma_2) \wedge (\beta_2 + \gamma_1)$.
\end{proposition}

\begin{remark}\label{rem:prod}
In general, for $F_i \in T_{<\gamma_i}$ one does not have $F_1 \ast F_2 \in T_{<\gamma}$
but only $F_1 \ast F_2 \in T_{<(\gamma_1+\gamma_2)}$. In order to avoid this, one can 
\textit{define} the product $\ast \colon \CD^{\gamma_1}\times \CD^{\gamma_1} \to \CD^{\gamma}$
as
\begin{equ}
(F_1 \ast F_2)(z) = \CQ_{<\gamma} \bigl(F_1(z) \ast F_2(z)\bigr)\;,
\end{equ}
where $\CQ_{<\gamma}$ is the projection onto $T_{<\gamma}$. It is a slightly lengthier exercise to verify
that Proposition~\ref{prop:prod} still holds with this modified definition, so we will use this from now on.
\end{remark}

Note here that $F_i$ can be seen as having two distinct ``regularities''. 
Its ``descriptive regularity'' $\gamma_i > 0$ measures how well $\CR F_i$ can be described in
terms of the underlying model, while its ``scaling regularity'' $\beta_i \le 0$ measures the 
behaviour of $\CR F_i$ around any given point in the sense that 
$\bigl(\CR F_i\bigr)(\phi_x^\lambda) \lesssim \lambda^{\beta_i}$. 
Note that the usual H\"older regularity
of a continuous function is its descriptive regularity (with respect to the polynomial model),
while its scaling regularity is always $0$. On the other hand, if we lift a distribution $\xi$ as above to a regularity structure via a
suitable symbol $\sXi$, then the function $F(x) \mapsto \sXi$ satisfying $\CR F = \xi$
has infinite descriptive regularity, while its scaling regularity coincides with the usual
(negative) H\"older regularity of $\xi$.
Proposition~\ref{prop:prod} can be interpreted as
a rigorous formulation of the intuition that ``multiplication by a distribution of regularity $-\alpha$
behaves like taking $\alpha$ derivatives''.

\begin{remark}
At the algebraic level, it is not difficult to see that given $V_1$ and $V_2$ it is always possible
to extend $T$ to some larger structure space $\tilde T$
 while keeping $G$ fixed in such a way that on has a product $\ast \colon V_1\times V_2 \to \tilde T$ 
and an action of $G$ on $\tilde T$ satisfying \eqref{e:algProd}. It is however \textit{not}
clear in general how to extend a given model on $T$ to all of $\tilde T$. 
It was shown in \cite[Prop.~4.11]{Regularity} that it is actually always \textit{possible} to do so, but
while the construction given there is explicit, it relies on several arbitrary choices making this
extension not canonical, except for homogeneous elements with degrees adding up to a strictly positive number. 
For example, in the situation of Remark~\ref{rem:notalg}, that construction
simply gives $\Pi_x \sXi_i \sXi_j = 0$, which is not very satisfactory at all.
In the particular case where the $\xi_i$ actually happen to be smooth, one would of course rather like
to set $\Pi_x \sXi_i \sXi_j = \xi_i\xi_j$.
\end{remark}

At this stage, our construction is already sufficiently rich to allow us to give a slightly cleaner 
formulation of the argument in Section~\ref{sec:Phi4}. Indeed, ignoring for a moment the fact that 
one would like to impose an initial condition at time $0$, we reformulate \eqref{e:Phi4} as an integral
equation:
\begin{equ}[e:integral]
\Phi = P \star \bigl(c\Phi - \Phi^3 + \xi\bigr) = P \star \bigl(c\Phi - \Phi^3) + \Psi\;.
\end{equ}
If we now consider a regularity structure built from the polynomial one by adding symbols $\sPsi^i$ for
$i \le 3$, we can then formulate this in a very natural way as a fixed point problem in $\CD^\gamma$:
\begin{equ}[e:Phi4abstr]
\Phi = \hat \CL \bigl(P \star \CR(c\Phi - \Phi^3)\bigr) + \sPsi\;,
\end{equ}
where we write $\hat \CL$ for the map turning a continuous $\CC^\gamma$ function into
an element of $\CD^\gamma$ as in \eqref{e:Taylor}.
Under what conditions can we even find a space $\CD^\gamma$ which is mapped into itself by the
right hand side of \eqref{e:Phi4abstr}? By our definition, any solution $\Phi$ will be of the form
$\Phi = \sPsi + R$ for some remainder $R$ taking values in the polynomial sector, so that 
$\Phi^3$ does at least have a natural meaning in our regularity structure. Furthermore, if $\deg \sPsi = -\alpha$,
and $\Phi \in \CD^\gamma$, we can apply Proposition~\ref{prop:prod} twice to find that 
$\Phi^3 \in \CD^{\gamma-2\alpha}$, so that we certainly want to be able to take $\gamma > 2\alpha$ for
the argument of $\CR$ to be sufficiently regular.

On the other hand, it follows from \eqref{e:boundModel} and Theorem~\ref{theo:reconstr} that one typically has 
$\CR \Phi^3 \in \CC^{-3\alpha}$ but no better so that, by the classical Schauder estimates, 
$\hat \CL P\star \CR \Phi^3 \in \CD^{2-3\alpha}$ but no better. These two conditions can be satisfied
simultaneously provided that $2\alpha < 2-3\alpha$, which restricts us to $\alpha < {2\over 5}$.
However, if $\xi$ is space-time white noise, then it only belongs to $\CC^{-\beta}$ for $\beta > {d+2\over 2}$,
thus enforcing $\alpha > {d-2\over 2}$ which restricts us to dimensions $d < {14\over 5}$.
In particular, this construction has no chance of covering the interesting case $d=3$. 

The problem is that, in \eqref{e:Phi4abstr}, we throw away all of the control that we have on
$\Phi^3$ the moment we apply the reconstruction operator $\CR$. Instead, we would like to
construct an operator $\CP$ which lifts the operation $P\star$ to the setting of modelled distributions 
in the sense that
\begin{equ}[e:intertwine]
\CR \CP F = P\star \CR F\;,
\end{equ}
but that furthermore satisfies a type of Schauder
estimate that increases both the descriptive and the scaling regularities by $2$
(modulo the limitation that the scaling regularity can never become strictly positive).

In order to build such an operator, we need to be able to describe the local behaviour
of $P\star \CR F$ up to order $\gamma + 2$, provided $F \in \CD^\gamma$. Again, this requires
our regularity structure to be sufficiently rich. It is natural to expect that if the local 
behaviour of $\CR F$ involves $\Pi_x \tau$ for some $\tau$, then that of $P\star \CR F$
will involve $P \star \Pi_x \tau$. Na\"\i vely, this would suggest that we should 
consider regularity structures endowed with a map $\CI$ such that 
\begin{equ}[e:degI]
\deg \CI \tau = \deg\tau + 2\;,
\end{equ}
and only consider models with the property that
\begin{equ}
\Pi_x \CI \tau = P \star \Pi_x \tau\;.
\end{equ}
Unfortunately, this is not compatible with the bound \eqref{e:boundModel}: even if $\deg \tau > -2$ for example,
there is no reason whatsoever why $P \star \Pi_x \tau$ should vanish at $x$, which is imposed by 
\eqref{e:boundModel}. Instead, the idea is to \textit{force} $\Pi_x \CI \tau$ to vanish at the right
order by considering models satisfying
\begin{equ}[e:admissible]
\Pi_x \CI \tau = P \star \Pi_x \tau - \sum_{|k| < \deg \CI\tau}{(\cdot -x)^k\over k!} \bigl(D^kP \star \Pi_x \tau\bigr)(x) \;.
\end{equ}
This immediately raises two questions:
\begin{enumerate}
\item Since $D^kP \star \Pi_x \tau$ is a genuine distribution in general, what does it mean to 
evaluate it at $x$?
\item Is it always possible to construct a regularity structure and model with these properties
and is this sufficient to construct the operator $\CP$?
\end{enumerate}
The first question can be addressed by realising that, by self-similarity, there exists a compactly supported test function $\phi$ so that the heat kernel $P$ can be written as
\begin{equ}
D^k P = \sum_{\log \lambda \in \Z} \lambda^{2 - |k|} (D^k\phi)_0^{\lambda}\;,
\end{equ} 
so that, at least formally, one has
\begin{equ}
\bigl(D^kP \star \Pi_x \tau\bigr)(x)
= \sum_{\log \lambda \in \Z} \lambda^{2 - |k|} \bigl(\Pi_x \tau\bigr)((D^k\phi)_x^{\lambda})\;.
\end{equ}
If we ignore the contributions coming from $\lambda > 1$ (one can deal with these separately or
replace the heat kernel by a truncated version), we see that this sum is guaranteed to converge precisely
when $\alpha + 2-|k| > 0$, i.e.\ when $|k| < \deg\CI\tau$.

Regarding the second question, it is shown in \cite[Sec.~4]{Regularity} that it does indeed
have a positive answer in the following sense. 

\begin{theorem}\label{theo:int}
Given a regularity structure $(T,G)$ extending the polynomial
structure $(\bar T, \bar G)$ and a subspace
$V \subset T$ (block-diagonal and invariant under $G$), one can further extend $(T,G)$ to a 
regularity structure $(\tilde T,\tilde G)$
on which there exists a map $\CI\colon V \to \tilde T$ satisfying \eqref{e:degI} and such that,
for every $\Gamma \in \tilde G$ and $\tau \in V$, one has
\begin{equ}[e:int]
\Gamma \CI \tau - \CI \Gamma \tau \in \bar T\;.
\end{equ}
Furthermore, every model on $(T,G)$ can be canonically and continuously
extended to a model on $(\tilde T, \tilde G)$ satisfying \eqref{e:admissible}.
Finally, writing $\CD^\gamma(V)$ for the space of $\CD^\gamma$ functions with values in $V$, 
one can construct a continuous operator $\CP \colon \CD^\gamma(V) \to \CD^{\gamma+2}$
such that \eqref{e:intertwine} holds and such that furthermore
\begin{equ}
\bigl(\CP F\bigr)(x) - \CI \bigl(F(x)\bigr) \in \bar T\;. 
\end{equ}
\end{theorem}

\subsection{Construction of regularity structures for SPDEs}
\label{sec:constr}

We now have all the ingredients in place to construct the space $\CX$ and the map $\hat \CM$ 
appearing in \eqref{e:wanted} and \eqref{e:defMA}. To construct the regularity structure, start
with symbols $\sone$, $\sX_i$ (one for every space-time coordinate) and $\sXi_j$ (one for every driving noise)
with degrees given as above. From these, we then build a collection $\CW$ of symbols by 
inductively postulating that
\begin{equs}[4]
\tau_1,\tau_2 &\in \CW \qquad&&\Rightarrow&\qquad \tau_1\tau_2 &\in \CW\;, \quad&\quad \deg \tau_1\tau_2 &= \deg\tau_1 + \deg\tau_2\\
\tau&\in \CW \qquad&&\Rightarrow&\qquad \CI(\tau) &\in \CW\;, \quad&\quad \deg \CI(\tau) &= 2+\deg\tau\;.
\end{equs}
We put on $\CW$ the equivalence relation $\sim$ generated by $\sone \tau \sim \tau$, 
$\tau_1\tau_2\sim \tau_2\tau_1$, $(\tau_1\tau_2)\tau_3 \sim \tau_1(\tau_2\tau_3)$, 
and by also postulating that $\tau_1 \sim \tau_2$ implies both
$\CI(\tau_1) \sim \CI(\tau_2)$ and $\tau_1\tau \sim \tau_2\tau$.
We then write $\tilde \CW = \CW / \sim$ so that $\tilde \CW$ is a graded commutative monoid endowed
with a map $\CI$.

The set $\tilde \CW$ is much too large to be a good set of basis vectors for a regularity structure 
since the corresponding set of degrees is typically neither discrete nor bounded from below.
One therefore selects the smallest subset $\CF \subset \tilde \CW$ which is sufficiently large so that
it is possible to formulate the SPDE under consideration as a fixed point problem in the regularity structure 
built by taking $T = \scal{\CF}$ and for $G$ the group of all linear operators on $T$ such that 
$\Gamma \sX_i - \sX_i \in \scal{\sone}$, $\Gamma \sXi_i = \sXi_i$, and 
satisfying \eqref{e:triang}, \eqref{e:algProd} and \eqref{e:int}. (The set $\CF$ needs to be chosen sufficiently
large so that any linear map $\Gamma$ on $\scal{\tilde \CW}$ satisfying all of the above properties maps
$\scal{\CF}$ to itself.) 

An equation is said to be locally subcritical if the corresponding set $\CF$ is such that,
for every $\gamma \in \R$, there are only finitely many elements $\tau \in \CF$ with $\deg\tau < \gamma$,
which is a prerequisite for $(T,G)$ as above to fulfil the axioms of a regularity structure.
Let us go back to our running example of the $\Phi^4$ model. In this case, ignoring again the effect
of initial conditions, we would like to rewrite \eqref{e:integral} as a fixed point problem
in some $\CD^\gamma$ space as follows:
\begin{equ}[e:PhiInt]
\Phi = \CP \bigl(c\Phi - \Phi^3 + \sXi\bigr)\;.
\end{equ} 
As a consequence of Theorem~\ref{theo:int}, this equation is of the form
\begin{equ}[e:genericForm]
\Phi = \CI \bigl(c\Phi - \Phi^3 + \sXi\bigr) + (\ldots)\;,
\end{equ}
where $(\ldots)$ denotes terms belonging to $\bar T$. Let us write $\CF_0$ for the collection of basis vectors 
required to describe $\Phi$ and $\CF_0$ for those required to describe $\Phi^3$. We then certainly want
$\sX^k \in \CF_0$ for all $k$ as well as $\CI(\sXi) \in \CF_0$. Furthermore, whenever
$\tau \in \CF_1$, we want to have $\CI(\tau) \in \CF_0$. Regarding $\CF_1$, whenever
$\tau_1, \tau_2, \tau_3 \in \CF_0$, we want to have $\tau_1\tau_2\tau_3 \in \CF_1$, so that 
$\Phi^3 \in \scal{\CF_1}$ for all $\Phi \in \scal{\CF_0}$.

To simplify notations, let us introduce graphical notations where we denote $\sXi$ by a dot and
$\CI$ by a line, and where multiplication is denoted by joining graphs by their roots.
For example, we write $\CI(\sXi) = \<1>$, $\CI(\sXi)^2 = \<2>$, $\CI(\CI(\sXi)^2)\CI(\sXi)^2 = \<22>$, etc. 
With this notation, we then have
\begin{equs}
\CF_0 &= \{\<1>, \<01>, \<02>, \<03>,\ldots\}\;,\\
\CF_1 &= \{\<1>,\<2>, \<3>, \<11>,\<21>,\<31>,\<22>,\<32>,\ldots \}\;.
\end{equs}
If we set $\deg\sXi = -\alpha$, then we have  
$\deg \<1> = 2-\alpha$, $\deg \<3> = 6-3\alpha$, $\deg \<22> = 10-4\alpha$, etc.
One can verify that the condition for subcriticality is that $\deg \<3> > \deg\sXi$, namely that 
$\sXi$ is the term of lowest degree appearing on the right hand side of \eqref{e:PhiInt}.
This imposes that $\alpha > -3$ so that, since space-time white noise has scaling exponent
$-{d+2\over 2}$ the critical dimension is indeed $4$, which is consistent with the heuristic presented in Section~\ref{sec:renormSPDE}.

It should be clear that this procedure is rather robust and allows us to associate a regularity structure to
any locally subcritical stochastic PDE. 
Furthermore, the regularity structures built with this procedure are endowed with a canonical lift $\CL \colon \CB_\infty \to \CX$
(recall that $\CX$ is the space of models for our regularity structure) by setting
\begin{equ}
\Pi_x \sXi_i = \xi_i\;,\qquad (\Pi_x \sX_i)(z) = z_i - x_i\;,
\end{equ}
and then extending this to all of $\CF$ by imposing that 
\begin{equ}[e:multcanon]
\Pi_x (\tau \bar \tau) = (\Pi_x \tau)\cdot(\Pi_x\bar \tau)\;,
\end{equ}
as well as \eqref{e:admissible}. There is also a natural (essentially unique) way of choosing $\Gamma_{xy} \in G$
such that \eqref{e:alg} holds and it is possible to verify that this choice of $(\Pi, \Gamma)$ satisfies the
bound \eqref{e:boundModel} as well as the required analytical bound on $\Gamma_{xy}$, provided that $\xi$ is smooth.

\section{Renormalisation}

Let us quickly summarise the situation so far. 
Given a class of stochastic PDEs, we have just seen how to build a regularity structure in which 
we can use the calculus developed in Section~\ref{sec:calculus} in order to reformulate the equation as a fixed
point problem in some space $\CD^\gamma$. In general, this equation may not have a solution but, if one
considers periodic (in space) situations, it was shown in \cite{Regularity} that it is possible to define weighted versions
of the spaces $\CD^\gamma$ in which the right hand side of \eqref{e:PhiInt} (or the analogous mild formulation 
of the SPDE under consideration, taking furthermore the effect of the initial condition into account) 
defines a contraction when restricted
to short enough time intervals, so that it admits a solution $\hat\CM$ which depends continuously on
the underlying model $Z \in \CX$.

Via the reconstruction operator $\CR$, $\hat\CM$ then defines a solution map $\CM_A$ by \eqref{e:defMA}. 
This construction implements the strategy outlined in Section~\ref{sec:proof}: the lift $\CL$ is given by the 
canonical lift as above
with the left inverse $\pi$ given trivially by $(\pi(\Pi,\Gamma))_i = \Pi_0 \sXi_i$. This is indeed surjective on 
$\CB$ provided that we choose $\deg\sXi_i = \alpha_i - \kappa$, since the bound \eqref{e:boundModel} is equivalent to
the $\CC^\alpha$ norm cases where $\Gamma_{xy}$ acts trivially on $\tau$ (as it does on $\sXi_i$ by construction). 
Furthermore, the identity \eqref{e:multcanon} combined with Remark~\ref{rem:pointwise} implies that, 
for smooth models, one has
\begin{equ}[e:canonical]
\CR(F\cdot G) = \CR F \cdot \CR G\;,
\end{equ}
which in turn, when combined with \eqref{e:intertwine}, shows 
that the first identity in \eqref{e:wanted} holds.

\subsection{Alternative representation of models}

The only remaining problem is that, as we have already seen in the calculations of Section~\ref{sec:Phi4bis}, 
the canonical lift $\xi \mapsto \CL(\xi)$ does not extend continuously to all $\xi \in \CB$. 
Instead, we would like to perform a suitable renormalisation procedure as described in Section~\ref{sec:proof}.
Before we do this, we note that the canonical lift $\CL$ can alternatively be described in a slightly different way.
Given a smooth $\xi \in \CB_\infty$, we obtain a representation $\PPi$ of our regularity structure by
setting similarly to above $(\PPi \one)(x) = 1$, $(\PPi \sX_i)(x) = x_i$, $(\PPi \sXi_i)(x) = \xi_i(x)$,
and then extending this inductively to all of $T$ by
\begin{equ}[e:canonLift]
\PPi (\tau \bar \tau) = (\PPi \tau)\cdot(\PPi \bar \tau)\;,\qquad \PPi \CI(\tau) = P\star \PPi \tau\;.
\end{equ}
This is very similar to the construction of $\CL(\xi)$ given above, except that we dropped the polynomial
correction term in \eqref{e:admissible}.
One then has the following alternative characterisation of $\CL(\xi)$ which is relatively straightforward
to show by induction.

\begin{proposition}\label{prop:canonical}
Given $\xi \in \CB_\infty$, define $\PPi$ as above. Then, for every $x$, there exists an essentially unique 
$F_x \in G$ such that $(D^k\PPi F_x \tau)(x) = 0$ for every $\tau \in \CF$ and every $|k| < \deg\tau$.
Furthermore, writing $\CL(\xi) = (\Pi,\Gamma)$, one has 
\begin{equ}[e:defModel]
\Pi_x = \PPi F_x\quad\text{and}\quad \Gamma_{xy} = F_x^{-1}F_y\;.
\end{equ}
\end{proposition}

\begin{remark}
By ``essentially unique'', we mean that $F_x$ is uniquely determined by $\PPi$ for generic 
choices of $\PPi$ (and therefore $\xi$). In special cases as when $\xi = 0$ this is of course
not the case. See \cite[Prop.~6.3]{BHZ} for a more precise formulation.
The condition $(D^k\PPi F_x \tau)(x) = 0$ for $|k| < \deg\tau$ combined with the smoothness
of $\PPi$ immediately implies \eqref{e:boundModel}.
\end{remark}

\begin{remark}\label{rem:renormPositive}
For arbitrary $\PPi \colon T \to \CC^\infty$, we can still find a (generically unique)
$F_x^{\PPi} \in G$ such that the condition $(D^k\PPi F_x \tau)(x) = 0$ for $|k| < \deg\tau$
holds for all $\tau$ of the form $\tau = \CI(\sigma)$ for some $\sigma \in \CF$.
This does however \textit{not} imply in general that this bound (and therefore \eqref{e:boundModel})
holds for arbitrary $\tau$!
\end{remark}

Write now $\CX_\infty$ for the collection of all smooth maps
\begin{equ}[e:admissibleBasic]
\PPi \colon T \to \CC^\infty\quad \text{such that}\quad \PPi \CI(\tau) = P\star \PPi \tau
\end{equ}
 and such that 
one can find a map $x\mapsto F_x \in G$ such that \eqref{e:defModel}
defines a model $(\Pi,\Gamma) \in \CX$. 
While Proposition~\ref{prop:canonical} states that the canonical representation $\PPi$ 
built from $\xi \in \CB_\infty$ does belong to $\CX_\infty$, 
it is \textit{not} the case in general that $\CX_\infty$ contains all 
maps $\PPi$ satisfying \eqref{e:admissibleBasic}.

Since every smooth model is built from a map $\PPi$ as in \eqref{e:defModel}
(take for example $\PPi = \Pi_0$ and $F_x = \Gamma_{0,x}$), we henceforth describe models
via such maps. 
Since the maps $\PPi$ are linear, it is natural in view of Section~\ref{sec:proof} to look
for a group $\RR$ of linear maps $M \colon T \to T$ acting on $\CX_\infty$ on the right by
$\PPi \mapsto \PPi \circ M$.
We furthermore want our renormalisation procedure to keep the meanings of the symbols 
$\sone$, $\sX$ and $\sXi$ intact and to  preserve the second identity in \eqref{e:canonLift}.
We therefore look for maps $M \colon T \to T$ satisfying the following properties
\begin{enumerate}[itemsep=0pt,parsep=0pt,topsep=0.2em]
\item One has $M \sone = \sone$ and $M \sXi_i = \sXi_i$ for all $i$.
\item One has $M \sX_i \tau = \sX_i M\tau$ and $M \CI(\tau) = \CI(M\tau)$ for all $\tau \in T$.
\item For every $\PPi \in \CX_\infty$, one has $\PPi \circ M \in \CX_\infty$.
\end{enumerate}
The first two properties are straightforward to verify, but it is much less clear \textit{a priori}
how the last one can be checked.
A group $\RR$ of linear maps satisfying all of these properties was exhibited in \cite[Thm~6.28]{BHZ},
 (it is called $\CG_-$ in that article) where it was also shown that $\RR$ is sufficiently large to be able to find 
elements $g^\xi \in \RR$ such that \eqref{e:center} holds. (Here, $\CX_\infty$ plays the role
of the linear space $\hat \CX$ alluded to in Remark~\ref{rem:Xhat} and $\hat \CL$ is given by the canonical
lift $\xi \mapsto \PPi$ described above.)

\subsection{Description of the renormalisation procedure}
\label{sec:RG}

We do not give a full description of the group $\RR$ here, but rather a slightly simplified one
corresponding to a subgroup of the group $\CG_-$ described in \cite{BHZ}. The idea
is that each element $M_g \in \RR$ is described by a map $g\colon \CF_- \to \R$, where
$\CF_- \subset \CF$ denotes the subset of $\CF$ consisting of symbols with strictly negative degree
and not containing any of the symbols $\sX_i$.
As in Section~\ref{sec:constr}, we identify elements of $\CF$ with trees $\tree = (V,E)$ whose nodes $V$
are endowed with labels describing the index $i$ of the corresponding term $\sXi_i$.
If the symbol in question also contains powers of $\sX$, then we encode these as additional 
labels on the nodes of the tree, so that nodes are indexed by elements in $\N^d \times \{0,\ldots,m\}$ where
a noise index of $0$ denotes the absence of noise.
For example, $\sX_2\CI(\sXi_1)\CI(\sXi_2 \sX_1)$ is identified with the tree $\<2>$, where the
two leaves have labels $((),1)$ and $((1),2)$ respectively, while the root has label $((2),0)$. 
(We denote multiindices by unordered tuples.)

Given such a labelled tree $\tree$, we write $\bar \tree \subset \tree$ for a subgraph
consisting of a collection of edges $\bar E \subset E$, as well as all the vertices $\bar V \subset V$
incident to $\bar E$. Given $\bar \tree \subset \tree$, we denote by $\iota(\bar \tree)$ the
element of the free algebra $\sscal{\CF_-}$ generated by $\CF_-$ by identifying each 
connected component of $\bar \tree$ with an element of $\CF_-$ (ignoring the values of the $\sX$-component
of the labels of the nodes of $\bar \tree$) and multiplying them in $\sscal{\CF_-}$. If one of the 
connected components of $\bar \tree$
happens to have positive degree, we set $\iota(\bar \tree) = 0$.
We naturally map the empty subgraph to the unit of $\sscal{\CF_-}$.

Given $\bar \tree \subset \tree$, we also denote by $\tree / \bar \tree$ the new element of $\CF$
obtained by contracting each connected component $A$ of $\bar \tree$ to a single vertex.
The label of this vertex is given by $(k,0)$, where $k$ is the sum of all the $\sX$-labels of 
vertices of $A$. For example, we have
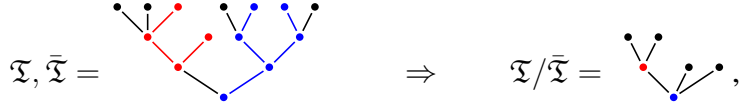
\begin{equ}
\tree, \bar \tree =   
\begin{tikzpicture}[scale=0.2,baseline=0.2cm]
        \node at (0,0)  [dot,blue] (root) {};
         \node at (-7,6)  [dot] (leftll) {};
          \node at (-5,4)  [dot,red] (leftl) {};
          \node at (-3,6)  [dot,red] (leftlr) {};
          \node at (-5,6)  [dot] (leftlc) {};
      \node at (-3,2)  [dot,color=red] (left) {};
         \node at (-1,4)  [dot,red] (leftr) {};
         \node at (1,4)  [dot,blue] (rightl) {};
          \node at (0,6)  [dot] (rightll) {};
           \node at (2,6)  [dot,blue] (rightlr) {};
           \node at (5,4)  [dot,blue] (rightr) {};
            \node at (4,6)  [dot,blue] (rightrl) {};
        \node at (3,2) [dot,blue] (right) {};
         \node at (6,6)  [dot] (rightrr) {};
        
        \draw[kernel1] (left) to node [sloped,below] {} (root);
        \draw[kernel1,red] (leftl) to
     node [sloped,below] {}     (left);
     \draw[kernel1,red] (leftlr) to
     node [sloped,below] {}     (leftl); 
     \draw[kernel1] (leftll) to
     node [sloped,below] {}     (leftl);
     \draw[kernel1,red] (leftr) to
     node [sloped,below] {}     (left);  
        \draw[kernel1,blue] (right) to
     node [sloped,below] {}     (root);
      \draw[kernel1] (leftlc) to
     node [sloped,below] {}     (leftl); 
     \draw[kernel1,blue] (rightr) to
     node [sloped,below] {}     (right);
     \draw[kernel1] (rightrr) to
     node [sloped,below] {}     (rightr);
     \draw[kernel1,blue] (rightrl) to
     node [sloped,below] {}     (rightr);
     \draw[kernel1,blue] (rightl) to
     node [sloped,below] {}     (right);
     \draw[kernel1,blue] (rightlr) to
     node [sloped,below] {}     (rightl);
     \draw[kernel1] (rightll) to
     node [sloped,below] {}     (rightl);
     \end{tikzpicture} 
     \qquad\Rightarrow\qquad \tree / \bar \tree =
   \begin{tikzpicture}[scale=0.2,baseline=0.2cm]
        \node at (0,0)  [dot,blue] (root) {};
          \node at (-2,2)  [dot,red] (leftl) {};
         \node at (-3,4)  [dot] (leftll) {};
          \node at (-1,4)  [dot] (leftlc) {};
           \node at (1,2)  [dot] (rightll) {};
         \node at (3,2)  [dot] (rightrr) {};
        
       \draw[kernel1] (leftl) to node [sloped,below] {} (root); ;
        \draw[kernel1] (leftll) to node [sloped,below] {} (leftl); ;
        \draw[kernel1] (leftlc) to
     node [sloped,below] {}     (leftl);
     \draw[kernel1] (rightll) to
     node [sloped,below] {}     (root); 
     \draw[kernel1] (rightrr) to
     node [sloped,below] {}     (root);
     \end{tikzpicture}\;,
\end{equ}
where $\bar \tree$ is the coloured subgraph of $\tree$, with the two connected components
drawn in blue and red respectively.
If the labels of the red vertices in $\tree$ are labelled by $(k_i,j_i)$, then the label
of the single red vertex in $\tree / \bar \tree$ is $(k,0)$ with $k = \sum_i k_i$.

Given $g \colon \CF_- \to \R$, we then also write $g$ for its unique extension to a character 
of $\sscal{\CF_-}$. The above definitions then allow us to define a map $M_g \colon T \to T$
by
\begin{equ}[e:actionAbstract]
M_g \tree = \sum_{\bar \tree \subset \tree} g(\iota(\bar \tree))\, \tree / \bar \tree\;.
\end{equ}
It is not very difficult to verify that one has $M_f \circ M_g = M_{f \circ g}$,
where the composition rule between maps $f ,g\colon \CF_- \to \R$ is given by
\begin{equ}[e:group]
(f\circ g)(\tree) =  \sum_{\bar \tree \subset \tree} g(\iota(\bar \tree))\, f(\tree / \bar \tree)\;,
\end{equ}
again with the conventions that $g$ is extended to all of $\sscal{\CF_-}$ as above and that 
$f(\tree / \bar \tree) = 0$ if its argument happens to be of positive degree.
Thanks to the fact that, as soon as $\bar \tree \neq \emptyset$,
$\tree / \bar \tree$ has strictly less edges than $\tree$, each operator $M_f$ 
differs from the identity by a (locally) nilpotent operator, so that the Neumann series
for its inverse converges. One can show that this inverse is of the form $M_{f^{-1}}$ for some
character $f^{-1}$, and that there exists a linear map $\CA\colon \scal{\CF_-} \to \scal{\CF_-}$
such that $f^{-1} = f\circ \CA$. (As a matter of fact, this endows $\scal{\CF_-}$ with a Hopf
algebra structure and \eqref{e:actionAbstract} turns $T$ into a left comodule for $\scal{\CF_-}$.)

\begin{remark}
Recall that the group $\RR$ constructed here is in general a proper subgroup of 
the full group $\CG_-$ described in \cite{BHZ} which is required for the renormalisation
procedure to work in general. The full group is constructed similarly to above, but without
the restriction that elements of $\CF_-$ have vanishing $\sX$-labels. In this case however,
formula \eqref{e:actionAbstract} has to be replaced by a more complicated version which
includes a suitable action on those labels. Since this would only distract from the gist
of the argument, we refrain from presenting this here.
\end{remark}

\subsection{Action on the space of models}

Denoting by $\RR$ the group of all characters of $\sscal{\CF_-}$ endowed with the composition rule 
\eqref{e:group}, we obtained in \cite[Thm~6.28]{BHZ} the following result.

\begin{theorem}
The formula $(f,\PPi) \mapsto \PPi \circ M_f$ yields a (right) group action
of $\RR$ onto $\CX_\infty$ which extends to a continuous group action on $\CX$.
\end{theorem}

The purpose of this section is to give an idea how the proof of this result works. 
Assume for a moment that the map $(g,F) \mapsto  F^g = M_g^{-1} F M_g$ yields a group action of
$\RR$ onto $G$ by a group automorphisms such that,
for every $\PPi \colon T \to \CC^\infty$ and every $g \in \RR$, writing $F_x^{\PPi}$ for the element of
$G$ described in Remark~\ref{rem:renormPositive}, one has 
\begin{equ}[e:wantedFx]
M_g F_x^{\PPi M_g} = F_x^{\PPi} M_g\;,
\end{equ}
or equivalently $(F_x^{\PPi})^g =  F_x^{\PPi M_g}$.

If this were the case then, writing $\PPi^g = \PPi M_g$, we would have the identity 
\begin{equ}[e:Pixg]
\Pi^g_x \eqdef \PPi^g F_x^{\PPi^g} = \PPi M_g (F_x^{\PPi})^g
 = \PPi F_x^{\PPi} M_g = \Pi_x M_g\;.
\end{equ}
Since, by the definition \eqref{e:actionAbstract}, $M_g \tau - \tau \in T_{> \deg \tau}$,
it immediately follows that $\PPi \in \CX_\infty$ implies $\PPi^g \in \CX_\infty$.
Unfortunately, it turns out that the map $(g,F) \mapsto  F^g = M_g^{-1} F M_g$ does
\textit{not} have these properties in general.

The way around this is to construct a larger regularity structure $(T^\ex, G^\ex)$ which
extends $(T,G)$ in the sense that $T \hookrightarrow T^\ex$, $G^\ex \twoheadrightarrow G$,
and the action of $G^\ex$ on $T^\ex$ leaves $T$ invariant and coincides with that of $G$ there.
This extension furthermore comes with a natural projection $\pi \colon T^\ex \to T$ as well
as with maps $\PPi^\ex \colon F_x^\ex \in G^\ex$ such that the collection of maps
$(\Pi_x^\ex,\Gamma_{xy}^\ex)$ given by
\begin{equ}
\Pi_x^\ex = \PPi^\ex F_x^\ex\;,\qquad \Gamma_{xy}^\ex = (F_x^\ex)^{-1}F_y^\ex\;,
\end{equ}
form a model for $(T^\ex, G^\ex)$ in the case where $\PPi^\ex$ is of the form
$\PPi^\ex = \PPi \pi$ for $\PPi$ the canonical lift of a smooth noise $\xi$.
The first step of the proof is then to show that the map $\iota\colon \PPi \mapsto \PPi \pi$ 
extends to a continuous injection of $\CX_\infty$ into $\CX_\infty^\ex$.
The second step is to show that $\RR$ also has a natural action $g \mapsto M_g^\ex$ 
on $T^\ex$ which is compatible with that on $T$ in the sense that 
\begin{equ}[e:actionProj]
\pi M_g^\ex = M_g\pi\;.
\end{equ}
In particular, this implies that $\iota( \PPi M_g) = \iota(\PPi)M_g^\ex$.
The third step then consists in verifying that for the extended regularity structure,
it is indeed the case that \eqref{e:wantedFx} holds, thus yielding a right 
action of $\RR$ onto $\CX_\infty^\ex$. Finally, one concludes the argument by 
showing that this action leaves $\CX_\infty$ (viewed as a subset of $\CX_\infty^\ex$
via $\iota$) invariant.

It remains to describe the extended regularity structure. The idea is to 
consider a set $\tilde \CW^\ex$ consisting of all labelled trees like $\tilde \CW$,
but this time with labels in $\N^d \times \{0,\ldots,m\} \times \R_-$, i.e.\ we
add an additional ``extended'' $\R_-$-valued label to each node. We write 
$\pi \colon \tilde \CW^\ex \to \tilde \CW$ for the map that simply forgets the extended label.
Writing $\beta_u$ for the extended label of the node $u$, we also define the degree
of $\tree \in \tilde \CW^\ex$ by setting
\begin{equ}
\deg\tree = \deg \pi\tree + \sum_{u \in V} \beta_u\;.
\end{equ}
Given such a labelled tree $\tree \in \tilde \CW$ as well as a subgraph $\bar \tree \subset \tree$
as above, we then define $\iota(\bar \tree) \in \sscal{\tilde \CW^\ex}$ just as above, but we do keep 
the values of the extended labels on the extracted subgraph.
We also define $\tree \dslash \bar \tree$ as in Section~\ref{sec:RG}, but 
with the extended label of any new vertex $u$ in $\tree \dslash \bar \tree$ corresponding to the contraction
of a connected component $A \subset \tree$ of $\bar\tree$ given by $\beta_{u} = \deg \iota(A)$. 
We then define $\CF^\ex \subset \tilde \CW^\ex$ as the set of all $\tree$ for which 
there exists $\hat \tree \in \CF$ and $\bar \tree \subset \hat \tree$ such that
$\tree = \hat \tree \dslash \bar \tree$.

The structure space for the extended regularity structure is then defined as $T^\ex = \scal{\CF^\ex}$. 
This space still admits a natural (partial) product obtained by joining trees by their roots
as above and setting the extended label at the root of the resulting tree to be the sum of the
extended labels at the roots of the two factors. Similarly, it admits a map $\CI$ in the same way 
as above, setting the extended label of the newly created root to $0$.
We can furthermore introduce maps $L_\alpha$ that act on a labelled tree $\tree$ by adding $\alpha$
to the extended label of the root.
We then define the group $G^\ex$ of all linear operators $\Gamma$ on $T^\ex$
such that $\Gamma \sX_i - \sX_i \in \scal{\sone}$, $\Gamma \sXi_i = \sXi_i$, and 
satisfying \eqref{e:triang}, \eqref{e:algProd} and \eqref{e:int}. We furthermore
restrict $G^\ex$ to elements such that $\Gamma L_\alpha(\tree) = L_\alpha(\Gamma \tree)$ whenever $\tree \in \CF^\ex$
is such that one also has $L_\alpha(\tree)\in \CF^\ex$.
It is clear by construction that any such map does act as some element of $G$ on the 
subspace $T \subset T^\ex$, so that the resulting regularity structure does indeed extend $(T,G)$. 

This allows us to build an ``extended'' renormalisation group $\RR^\ex$ exactly as above, but keeping
track of the extended labels. In other words, we define $\CF_-^\ex \subset \CF^\ex$ as the subset of
all elements of $\CF^\ex$ of strictly negative degree and with vanishing $\sX$-labels (but possibly
non-vanishing extended labels). Elements of $\RR^\ex$ are then described by maps
$g \colon \CF_-^\ex \to \R$ and their action on $T^\ex$ is given as above by
\begin{equ}[e:actionAbstractex]
M_g \tree = \sum_{\bar \tree \subset \tree} g(\iota(\bar \tree))\, \tree \dslash \bar \tree\;.
\end{equ}
One can verify that the map $g \mapsto g^\ex = g \circ \pi$ 
allows to view $\RR$ as a subgroup of $\RR^\ex$, thus yielding an action of
$\RR$ on $T^\ex$ which does indeed satisfy \eqref{e:actionProj}.

We claim that this action now \textit{does} satisfy the identity \eqref{e:wantedFx}.
The main reason for this is that, by the definition of the operation $\dslash$ and the degree
on $\CF^\ex$, all terms appearing in the right hand side of \eqref{e:actionAbstractex}
have the same degree. As a consequence, one can easily verify that if we define
$\Pi_x^g = \Pi_x M_g^\ex$ as in \eqref{e:Pixg}, then it \textit{does} satisfy 
\eqref{e:admissible}. Since the latter determines the characters $F_x^{\PPi}$, one can 
use this to go backwards and show that \eqref{e:wantedFx} does indeed hold.

\subsection{Renormalised equations}

In this last section, we show how one constructs the `dual' right 
action of $\RR$ onto the space of equations in such a way that the identity \eqref{e:renEqu}
holds. Formalising the space $\CN$ of ``possible right hand sides'' is rather heavy notationally, so
we restrict ourselves again to the specific example of the $\Phi^4$ model, although the
argument is quite general. 

We consider the case of $d = 3$, so that we can choose $\alpha = {5\over 2} + \kappa$ for 
some sufficiently small $\kappa > 0$. In particular, one has $\deg\<1> = -{1\over 2}-\kappa$,
so that Proposition~\ref{prop:prod} tells us that we should look for solutions to \eqref{e:PhiInt}
in $\CD^\gamma$ for some $\gamma > 1 +2\kappa$. Choosing $\gamma = 1+3\kappa$ and $\kappa$ small
enough, we then note that as a consequence of \eqref{e:genericForm}, every solution to
the fixed point problem \eqref{e:PhiInt} is of the form
\begin{equ}[e:sol]
\Phi = \<1> + \phi \sone - \<03> - 3\phi \<02> + \scal{\nabla \phi,\sX}\;,
\end{equ}
for some continuous functions $\phi$ and $\nabla \phi$. 

\begin{remark}
By Remark~\ref{rem:prod}, the identity \eqref{e:sol} is \textit{exact}, not just some approximate identity 
``up to higher order terms''. The continuous functions $\phi$ and $\nabla \phi$ however are
unknown and to be determined by the fixed point problem. Note also that $\nabla \phi$
is \textit{not} the gradient of $\phi$, but can be interpreted as a kind of ``renormalised gradient''.
In fact, it follows from \eqref{e:Dgamma} and \eqref{e:admissible} that for every 
$\Phi \in \CD^\gamma$ of the form \eqref{e:sol},
one has the identity
\begin{equ}
(\nabla_i \phi)(z) = \lim_{h\to 0} h^{-1} \bigl(\phi(z+ he_i) - \phi(z) - \bigl(\Pi_z \<03>\bigr)(z+he_i)  - 3\phi(z)\bigl(\Pi_z \<02>\bigr)(z+he_i) \bigr)\;,
\end{equ}
with $e_i$ the $i$th unit vector in $\R^3$. The mere fact that there even \textit{exist} functions
such that this limit is finite for every $z$ is not obvious!
\end{remark}

As a consequence of \eqref{e:sol} and Remark~\ref{rem:prod}, one has the identity
\begin{equ}[e:Phi3]
\Phi^3 = \<3> + 3\phi \<2> + 3\phi^2\<1> + \phi^3\sone - 6\phi \<31> -3\<32> - 9\phi \<22> + 3\scal{\nabla \phi,\sX}\<2>\;.
\end{equ}
Again, this is true whatever the underlying model $\PPi$, as long as it is admissible in the sense that 
the second identity in \eqref{e:canonLift} holds, whether it is obtained from the canonical
lift of a smooth noise, or not. Of course the actual function / distribution $\CR \Phi$ represented by this
identity depends very much on $\PPi$, first because the reconstruction operator $\CR$ depends on
it and second because the function  $\phi$ (and therefore also $\nabla \phi$) depends on it.

Consider now the particular case when $\PPi$ is the canonical lift of some smooth noise $\xi$. 
In this case, it follows from \eqref{e:canonical} that $\CR(\Phi^3) = (\CR \Phi)^3$.
Let now $g$ be an element of $\RR$, write $\hat \PPi = \PPi \circ M_g$, and write $\hat \CR$ for the 
corresponding reconstruction operator. Making use of the extended model described in the previous
subsection, combined with Remark~\ref{rem:pointwise}, we can write this renormalised reconstruction operator as
\begin{equ}
\bigl(\hat \CR U\bigr)(z) = \bigl(\Pi_z^\ex M_g^\ex U(z)\bigr)(z) = \bigl(\CR M_g^\ex U\bigr)(z)\;.
\end{equ}
At this point, we choose $g\colon \CF_- \to \R$ as
\begin{equ}[e:gsimple]
g(\<2>) = -c_1\;,\qquad g(\<22>) = -c_2\;,
\end{equ}
and $g(\tau) = 0$ otherwise. It then follows from \eqref{e:actionAbstract} that one has
for example
\begin{equ}
M_g \<3> = \<3> - 3c_1 \<1>\;,\qquad
M_g \<31> = \<31> - 3c_1 \<11>\;,\qquad
M_g \<32> = \<32> - 3c_1 \<12> - c_1\<30> + 3c_1^2 \<01> - 3c_2\<1>\;.
\end{equ}
The action of $M_g^\ex$ is essentially the same, but with extended decorations with value $-1-2\kappa$
added on those vertices obtained by contracting an instance of $\<2>$.

Inserting this into \eqref{e:Phi3}, it follows that one has the identities
\begin{equ}
M_g^\ex \Phi \sim \Phi\;,\qquad M_g^\ex \Phi^3 \sim \Phi^3 - (3c_1 - 9c_2)\Phi\;,
\end{equ}
where we write $\tau \sim \bar \tau$ if $(\Pi_z^\ex\tau)(z)=(\Pi_z^\ex\bar\tau)(z)$.
We conclude that in this particular example one has 
\begin{equ}[e:renormPhi3]
\hat \CR(\Phi^3) = (\hat \CR \Phi)^3  - (3c_1 - 9c_2)\hat \CR \Phi\;.
\end{equ}
Take now as our space $\CN$ of nonlinearities all nonlinearities of the 
type $F_c(\Phi) =  c\Phi - \Phi^3$ for $c \in \R$.
It then follows from \eqref{e:renormPhi3}
 that, at least for the subgroup of $\RR$ given by elements of the type \eqref{e:gsimple} 
(which as a group is simply $(\R^2,+)$), one does indeed have the announced identity \eqref{e:renEqu} with
\begin{equ}
F_c g = F_{c + 3c_1 - 9c_2}\;.
\end{equ}

\begin{remark}
In the case $c_2 = 0$ (which is the case in dimension $2$), we see that we recover the 
calculation performed in \eqref{e:actionSmooth}. 
\end{remark}

\begin{remark}
The observant reader may have noticed that we have
\begin{equ}
g(\Phi^3) = -3c_1 \phi + 9c_2 \phi\;,
\end{equ}
which appears coincidentally to be ``essentially the same'' as the counterterm
$-3c_1\,\hat \CR \Phi + 9 c_2 \,\hat \CR \Phi$ appearing in \eqref{e:renormPhi3}.
This is actually \textit{not} a coincidence but holds in very wide generality,
see \cite{BCCH}. (Compare Eq.~2.15 with Def.~3.15 and Lemma~4.5, noting that the 
symmetry factor $S(\tau)$ is generated by the scalar product appearing in the left hand side
of Def.~3.15.)
This gives us a very simple way of deriving the renormalisation
procedure for any given class of stochastic PDEs.
\end{remark}

\bibliographystyle{Martin}

\bibliography{Takagi}

\end{document}